\documentclass[review,onefignum,onetabnum]{siamart171218}
\usepackage{tikz}
\hypersetup{colorlinks=true,
citecolor=blue,linkcolor=blue
}

\newcommand{\ad}{\mbox{ad}}

\usepackage{lipsum}
\usepackage{amsfonts}
\usepackage{graphicx}
\usepackage{epstopdf}
\usepackage{algorithmic}
\ifpdf
  \DeclareGraphicsExtensions{.eps,.pdf,.png,.jpg}
\else
  \DeclareGraphicsExtensions{.eps}
\fi

\newsiamremark{example}{Example}
\newsiamremark{remark}{Remark}
\newsiamremark{hypothesis}{Hypothesis}
\crefname{hypothesis}{Hypothesis}{Hypotheses}
\newsiamthm{claim}{Claim}

\headers{Variational integrators for optimal control of foldable drones}{L. Colombo, J. Giribet, and D. Mart\'in de Diego}

\title{Variational integrators on Lie groups for  optimal control of foldable drones\thanks{Submitted to the editors DATE.
\funding{L. Colombo and D. Mart{\'\i}n de Diego acknowledge financial support from the Spanish Ministry of Science and
Innovation, under grants  PID2022-137909NB-C21, TED2021-129455B-I00 and  CEX2023-001347-S funded by MCIN/AEI\-/10.13039\-/501100011033 and the LINC Global project from CSIC ``Wildlife Monitoring Bots'' INCGL20022. J. Giribet was partially supported by PICT-2019-2371 and PICT-2019-0373 projects from Agencia Nacional de Investigaciones Cient\'ificas y Tecnol\'ogicas, Argentina.}}}

\author{Leonardo J. Colombo\thanks{Centro de Automática y Robótica (CSIC-UPM), 
Carretera de Campo Real, km 0, 200, 28500 Arganda del Rey, Spain.
  (\email{leonardo.colombo@csic.es}).}
\and Juan I. Giribet\thanks{Universidad de San Andr\'es (UdeSA) and CONICET, Argentina. Calle Vito Dumas 284, B1644BID, Victoria, Buenos Aires, Argentina (\email{jgiribet@conicet.gov.ar}).}
\and David Mart\'in de Diego\thanks{Instituto de Ciencias Matematicas (CSIC), Calle Nicol\'as Cabrera 13-15, Cantoblanco, 28048, Madrid, Spain. (\email{david.martin@icmat.es}).}}

\usepackage{amsopn}

\makeatletter
\newcommand*{\addFileDependency}[1]{
  \typeout{(#1)}
  \@addtofilelist{#1}
  \IfFileExists{#1}{}{\typeout{No file #1.}}
}
\makeatother

\ifpdf
\hypersetup{
  pdftitle={Variational integrators on Lie groups for  optimal control of foldable drones},
  pdfauthor={L. J. Colombo, J. I. Giribet, and D. Mart\'in de Diego}
}
\fi

\begin{document}

\maketitle

\begin{abstract}
 Numerical methods that preserve geometric invariants of the system, such as energy, momentum or the symplectic form, are called geometric integrators. These include variational integrators as an important subclass of geometric integrators. The general idea for those variational integrators is to discretize Hamilton’s principle rather than the equations of motion and as a consequence these methods preserves some of the invariants of the original system (symplecticity, symmetry, good behavior of energy...). In this paper, we  construct variational integrators for control-dependent Lagrangian systems on Lie groups.  These integrators are derived via a discrete-time variational principle for discrete-time control-dependent reduced Lagrangians. We employ the variational integrator into optimal control problems for path planning of foldable unmanned aerial vehicles (UAVs). Simulation are shown to validate the performance of the geometric integrator.
\end{abstract}



\section{Introduction}
Since the emergence of computational methods, fundamental properties such as accuracy, stability, convergence, and computational efficiency have been considered crucial for deciding the utility of a numerical algorithm. Geometric numerical integrators are concerned with numerical algorithms that preserve the system's fundamental physics by keeping the geometric properties of the dynamical system under study. The key idea of the structure-preserving approach is to treat the numerical method as a discrete dynamical system which approximates the continuous-time flow of the governing continuous-time differential equation, instead of focusing on the numerical approximation of a single trajectory. Such an approach allows a better understanding of the invariants and qualitative properties of the numerical method. Using ideas from differential geometry, structure-preserving integrators have produced a variety of numerical methods for simulating systems described by ordinary differential equations preserving its qualitative features. In particular, numerical methods based on discrete variational principles \cite{mawest, Lall1} may exhibit superior numerical stability and structure-preserving capabilities than traditional integration schemes for ordinary differential equations.

Variational integrators are geometric numerical methods derived from the discretization of variational principles \cite{mawest, Lall1,  Hair}. These integrators retain some of the main geometric properties of the continuous systems, such as preservation of the manifold structure at each step of the algorithm, symplecticity, momentum conservation (as long as the symmetry survives the discretization procedure), and a good behavior of the energy function associated to the system for long time simulation steps. This class of numerical methods has been applied to a wide range of problems in optimal control \cite{ober2011discrete, leo,fer}, constrained systems \cite{Ley}, formation control of multi-agent systems \cite{colombo2021forced}, nonholonomic systems \cite{Cort,ferraro2008momentum}, accelerated optimization \cite{campos2021discrete}, non-autonomous systems \cite{colombo2023variational}, flocking control \cite{colombo2020forced} and motion planning for underactuated robots \cite{KM}, among many others. 

In this paper, we construct variational integrators for attitude control of foldable multirotors unmanned aerial vehicles (UAVs). One of the main advantages of multi-rotors over other UAVs such as fixed-wing UAVs, is their maneuverability. This makes them a more suitable platform to perform certain tasks, such as the inspection of civil structures, or navigating in narrow environments such as landslide areas. Motivated by different applications, efforts have been made to increase the maneuverability of these vehicles, emerging new designs capable of being reconfigured during flight \cite{derrouaoui2022comprehensive}. In \cite{falanga2018foldable}, \cite{hu2021design}, some of these reconfigurable vehicle designs have been experimentally validated, showing how it is possible to adapt their shape to access places where it would be impossible otherwise. In \cite{bucki2019design}, a multirotor design is proposed where retractable arms allow the vehicle to be launched from a platform, to then convert to a different flight configuration. These novel designs, with reconfigurable structures, present interesting challenges from the control design perspective.

In the designs proposed in \cite{falanga2018foldable}, \cite{papadimitriou2021geometry}, the vehicle adapts its shape according to the mission to be executed. For instance, if the vehicle needs to retract its arms to pass through a narrow area. Thus, the change in the vehicle's moment of inertia becomes a parameter, based on which the control law is based. Similarly, in \cite{Pose22}, a vehicle is proposed that adjusts its moment of inertia to compensate for a rotor failure. Similar to the previous studies, the alteration of the moment of inertia is taken as a parameter, and the control law adapts to this parameter variation.

In \cite{falanga2018foldable}, an adaptive LQR strategy is adopted for each possible configuration of the vehicle, which requires great computational power to obtain the control commands needed to stabilize the vehicle. In order to consider the constraints from the point of view of the control, in \cite{papadimitriou2021geometry} a model predictive control is proposed, and \cite{hu2021design} proposes a novel RL-based strategy for attitude stabilization. This strategy of adjusting the control law while considering the change in moment of inertia is also used in \cite{Zhai18}, where a system for on-orbit refueling is proposed to extend the lifespan of satellites lacking propellant. Throughout the propellant transfer, variations in inertia tensors occur, and an attitude control system is suggested, tailored to these moments of inertia changes. Unlike previous works, this paper treats the alteration of the UAV's moment of inertia not as a parameter but as a control action, incorporating moment of inertia changes within the objective function to minimize. Thus, the change in shape of the UAV becomes an additional control action. However, when viewed from a geometric standpoint, it produces distinct effects compared to external torque control inputs.

In this paper, we construct variational integrators for control-dependent Lagrangian systems and apply them to optimal control problems of foldable multirotor UAVs. In particular, in Section  \ref{sec:lg-background} we introduce the dynamics of reduced mechanical system on Lie groups and the Euler-Poincar\'e and Lie-Poisson equations. Further, in Section \ref{sec3.1} we study the dynamics of control-dependent Lagrangian systems, in particular, rigid bodies with controlled inertia, and we obtain the necessary conditions for optimality in an optimal control problem from an intrinsic geometric point of view. In Section \ref{sec: foldable} we present the model for the attitude dynamics of a foldable quadrotor considering the impact of the varying nature of the vehicle’s inertia and the associated optimal control problem for trajectory tracking. In Section \ref{sec: varint} we introduce variational integrators on Lie groups and in Section \ref{sec4.2} we construct variational integrators for optimal control problems of control-dependent Lagrangian systems. Finally in Section \ref{sec5} we derive variational integrators for the trajectory tracking of foldable multirotor UAVs and we show some simulation results. 

 \section{Dynamics of mechanical systems on Lie groups}\label{sec:lg-background}

A \textit{Lie group} is a finite dimensional smooth manifold $G$ that is a group and for which the operations of multiplication $(g,h)\mapsto gh$ for $g,h\in G$, and inversion, $g\mapsto g^{-1}$, are smooth.


Let $G$ be a Lie group with identity element $e\in G$. A \textit{left-action} of $G$ on a smooth manifold $Q$ is a smooth mapping $\Phi:G\times Q\to Q$ such that $\Phi(e,q)=q$, $\forall q\in Q$; $\Phi(g,\Phi(h,q))=\Phi(gh,q)$, $\forall g,h\in G, q\in Q$; and for every $g\in G$, $\Phi_g:Q\to Q$ defined by $\Phi_{g}(q):=\Phi(g,q)$ is a diffeomorphism. $\Phi:G\times Q\to Q$ is a \textit{right-action} if it satisfies the same conditions as for a left action
except that $\Phi(g,\Phi(h,q))=\Phi(hg,q)$ $\forall g,h\in G, q\in Q$. We often use the notation $\Phi_{g}(q)=\Phi(g,q)$ and say that $g$ acts on $q$. 

Let $G$ be a finite dimensional Lie group and let $\mathfrak{g}$ be its Lie algebra, $\mathfrak{g}:=T_{e}G$, the tangent space at the identity $e\in G$.  Let $\mathcal{L}_{g}:G\to G$ be the left translation of the element $g\in G$ given by $\mathcal{L}_{g}(h)=gh$ for $h\in G$. Similarly, ${\mathcal R}_g$ denotes the right translation of the element $g\in G$ given by ${\mathcal R}_{g}(h)=hg$ for $h\in G$. 
Their tangent maps  (i.e, the linearization or tangent lift) are denoted by $T_{h}\mathcal{L}_{g}:T_{h}G\to T_{gh}G$ and $T_{h}{\mathcal R}_{g}:T_{h}G\to T_{hg}G$, respectively. Similarly, the cotangent maps (cotangent lift) are denoted by $T_{h}^{*}\mathcal{L}_{g}: T^{*}_{h}G\to T^{*}_{gh}G$ and $T_{h}^{*}{\mathcal R}_{g}:T^{*}_{h}G\to T^{*}_{hg}G$, respectively. It is well known that the tangent and cotangent lifts are Lie group actions (see \cite{HolSch09}, Chapter $6$).

Let $\Phi_g:Q\to Q$ for any $g\in G$ a left action on $G$; a function $f:Q\to\mathbb{R}$ is said to be \textit{invariant} under the action $\Phi_g$, if $f\circ\Phi_g=f$, for any $g\in G$ (that is, $\Phi_g$ is a symmetry of $f$). 

Let $L\colon TG\to\mathbb{R}$ be a Lagrangian function. If we assume that the Lagrangian $L$ is $G$-invariant under the tangent lift of left translations, that is $L\circ T_{g}\mathcal{L}_{g^{-1}}=L$ for all $g\in G$, then it is possible to obtain a reduced Lagrangian $\ell\colon\mathfrak{g}\to\mathbb{R}$, where $$\ell(\omega) = L(g^{-1}g,T_{g}\mathcal{L}_{g^{-1}}(\dot{g}))= L(e,\omega).$$The reduced Euler--Lagrange equations, that is, the \textit{Euler--Poincar{\'e} equations} (see, e.g.,  \cite{HolSch09}, \cite{Blo03}), are given by the system of $n$ first-order ODE's  ($\dim G=n$)\begin{align}
\frac{d}{dt}\frac{\partial\ell}{\partial\omega} = \ad^{*}_{\omega}\frac{\partial\ell}{\partial\omega},\label{eq_ep_intro}
\end{align} where $\ad^{*}:\mathfrak{g}\times\mathfrak{g}^{*}\to\mathfrak{g}^{*}$, $(\omega,\Pi)\mapsto\ad^{*}_{\omega}\Pi$ is the \textit{co-adjoint operator} defined by $\langle\ad_{\omega}^{*}\Pi,\eta\rangle=\langle\Pi,\ad_{\omega}\eta\rangle$ for all $\eta\in\mathfrak{g}$ with $\ad:\mathfrak{g}\times\mathfrak{g}\to\mathfrak{g}$ the \textit{adjoint operator} given by $\ad_{\omega}\eta:=[\omega,\eta]$, where $[\cdot,\cdot]$ denotes the Lie bracket of vector fields on the Lie algebra $\mathfrak{g}$, and where $\langle\cdot,\cdot\rangle:\mathfrak{g}^{*}\times\mathfrak{g}\to\mathbb{R}$ denotes the so-called \textit{natural pairing} between vectors and co-vectors. For matrix Lie groups it is defined by $\langle\alpha,\beta\rangle:=\alpha\cdot\beta$ for $\alpha\in\mathfrak{g}^{*}$ and $\beta\in\mathfrak{g}$ where $\alpha$ is understood as a row vector and $\beta$ a column vector
(see \cite{HolSch09}, Section $2.3$ for details).

For the modelization of the attitude dynamics of a multi-rotor unmanned aerial vehicle (UAV) we need to add external torques to the Euler-Poincar\' e equations. For that, denote by $\mathfrak{so}(3)$ the Lie algebra of the special orthogonal Lie group $SO(3)$, that is, the set of $3\times 3$ skew-symmetric matrices, and consider the Lagrangian system $\ell: {\mathbb R}^3\equiv {\mathfrak {so}}(3)\rightarrow {\mathbb R}$ given by $\ell(w)=\frac{1}{2}I\omega\cdot\omega$ subject to external torques $F: G\times \mathfrak{so}(3) \to\mathfrak{so}(3)^ {*}$. Therefore, the equations are now   \begin{equation}\label{eq1}
\dot{\omega}=I^{-1}(F-\omega\times I\omega),
\end{equation} where $I$ is the inertia matrix considered as  possitive definite symmetric bilinear form on ${\mathfrak g}$ (or also considered as a linear map from ${\mathfrak g}$ to  ${\mathfrak g}^*$) and where $F: G\times \mathfrak{so}(3)\rightarrow\mathfrak{so}(3)^ {*}$ is an external torque. 

The Euler--Poincar{\'e} equations on $\mathfrak{g}$ \eqref{eq_ep_intro} together with the reconstruction equation $\omega = T_{g}L_{g^{-1}}(\dot{g})$ are equivalent to the Euler--Lagrange equations on $TG$. By assuming that the Legendre transformation for the reduced Lagrangian $\ell$ is a global diffeomorphism (i.e. $\ell$ is hyper-regular), then one can move to the Hamiltonian side. For that we define the reduced Hamiltonian $h\colon \mathfrak{g}^{*}\to\mathbb{R}$ given by
$
h(\Pi) = \langle\Pi,\omega(\Pi)\rangle-\ell(\omega(\Pi)),\nonumber
$ by using the reduced Legendre transformation $F\ell:\mathfrak{g}\to\mathfrak{g}^{*}$ given by $\langle F\ell(\omega),\eta\rangle=\langle\frac{\partial\ell}{\partial\omega}, \eta\rangle$ to express $\omega$ as a function of $\Pi$ by the implicit function theorem.  The Euler--Poincar{\'e} equations \eqref{eq_ep_intro} can then be written as the \textit{Lie--Poisson equations} on $\mathfrak{g}^{*}$, which are given by \begin{equation}\label{eq: lie-poisson}\dot{\Pi} = \ad^{*}_{\frac{\partial h}{\partial\Pi}}\Pi.\end{equation}


Using the fact that the angular momentum in the body frame can be written as $\Pi=I\omega$, equation \eqref{eq1} can be written as the Lie-Poisson equations subject to external torques \begin{equation}\label{controlled-Lie-Poisson}
    \dot{\Pi}=\bar{F}-I^{-1}\Pi\times\Pi.
\end{equation}
where  $\bar{F}(g, \Pi)\equiv F(g, I^{-1}\Pi)$.  

Note that equations \eqref{eq1} and \eqref{controlled-Lie-Poisson} are the (forced) Euler-Poincar\'e and  Lie-Poisson equations, respectively, for the reduced Lagrangian and Hamiltonian functions $\ell:\mathfrak{so}(3)\simeq\mathbb{R}^{3}\to\mathbb{R}$ and $h:\mathfrak{so}(3)^{*}\simeq\mathbb{R}^{3}\to\mathbb{R}$, given by $\ell(\omega)=\frac{1}{2}\langle I\omega,\omega\rangle$ and $h(\Pi)=\frac{1}{2}\langle\Pi,I^{-1}\Pi\rangle$, respectively.

\section{Optimal control for controlled inertia rigid bodies}\label{sec:ep-dynamics}

Inspired by the recent results of foldable drones \cite{falanga2018foldable} we aim to study rigid body attitude dynamics for systems with control-dependent inertia. 

\subsection{Necessary conditions for optimallity of controlled inertia rigid bodies}\label{sec3.1}
Next, we assume that the Lagrangian function is given now by $L:\mathfrak{so}(3)\times U\to\mathbb{R}$ as $$L(\omega,u)=\frac{1}{2}\langle I(u)\omega,\omega\rangle=\frac{1}{2}\omega^{T}I(u)\omega$$ where $U$ denotes the set of admissible (internal) controls. 



For simplified notation, we will denote $\mathfrak{so}(3)={\mathfrak g}$ in the sequel.  By Lagrange-d'Alembert principle \cite{Blo03}, the \textit{forced controlled Euler-Poincar\'e equations} for a Lagrangian $L:\mathfrak{g}\times U\to\mathbb{R}$ with external force $F:G\times\mathfrak{g}\times U\times\bar{U}\to \mathfrak{g}^{*}$ (see Appendix A for details) are given by \begin{equation}\label{eqforL}
I(u)\dot{\omega}+\left(\frac{\partial I}{\partial u}\dot{u}\right)\omega=F-\omega\times I(u)\omega,
\end{equation}
where $F=T_{e}^{*}L_gf(g,\omega,u,\tau)$ and 
where $\bar{U}$ is a set of admissible (external) controllers denoted by $\tau$. In addition, by defining the body angular momentum $\Pi=I(u)\omega$ we can write equations \eqref{eqforL} as the \textit{forced controlled Lie-Poisson equations} given by (see Appendix A for details) \begin{equation}\label{eq8}
\dot{\Pi}=\bar{F}-I^{-1}(u)\Pi\times\Pi.
\end{equation}

In the same manner than equations \eqref{eqforL}, equations \eqref{eq8} can be obtained by a variational principle for a reduced controlled Hamiltonian $h:\mathfrak{g}^{*}\times U\to\mathbb{R}$ and now the external force is $\bar{F}:G\times\mathfrak{g}^*\times U\times\bar{U}\to \mathfrak{g}^{*}$ (see Appendix A for details).

\subsubsection{Optimal control problem}\label{ocp} Along this paper we are interested on control an UAV to reach a desired equilibrium state $(g_d,\Pi_d)\in G\times\mathfrak{g}^{*}$. For this purpose we need to consider an optimal control problem with fixed final time $T>0$, $$\min\int_{0}^{T}C(g,\Pi,u,\dot{u},\tau)dt+c\Phi(g(T),\Pi(T)),\quad C:G\times\mathfrak{g}^{*}\times TU\times \bar{U}\to\mathbb{R}$$  subject to fixed boundary conditions $g(0)=g_0$, $g(T)=g_T$, $\Pi(0)=\Pi_0$, $\Pi(T)=\Pi_{T}$, $u(0)=u_0$, $u(T)=u_{T}$, together with the dynamic and kinematic equations \begin{align}\dot{\Pi}&=\bar{F}-I^{-1}(u)\Pi\times\Pi,\label{kinematics1}\\ \dot{g}&=gI^{-1}(u)\Pi.\label{kinematics2}\end{align} Here $\Phi:G\times\mathfrak{g}^{*}\to\mathbb{R}$ is a terminal cost (Mayer
term) and $c>0$ is a weight for the terminal cost. We assume that $\Phi$ is continuously differentiable function. 

\subsubsection{Necessary conditions for optimality} In the following, for simplicity in the computations, we will assume $I(u)$ is diagonal, $I(u)=\hbox{diag}(I_1(u),I_2(u),I_3(u))$. The optimality conditions for this problem are derived as follows. First consider the phase space $M:=(G\times\mathfrak{g}^{*})\times T^*{\mathfrak g}^*\times(TU\oplus T^{*}U)\times\bar{U}$ with local coordinates $(g,p_{\xi},\Pi,p_{\Pi},u,\dot{u},p_u,\tau)$ and  the Pontryagin Hamiltonian $\bar{H}_{M}:M\to\mathbb{R}$ given by \begin{equation}\label{opcHam}
\bar{H}_M=C(g,\Pi,u,\dot{u},\tau)+\langle p_{\Pi}, \bar{F}-I^{-1}(u)\Pi\times\Pi\rangle 
+\langle p_{\xi}, I^{-1}(u)\Pi\rangle
\end{equation} 
To derive Hamilton equations for $\bar{H}_M$ we consider the (presymplectic) closed $2$-form on $M$ given by $\Omega_M=\Omega_{G\times\mathfrak{g}^*}+\Omega_{{\mathfrak g}^*}+\Omega_U$ where $\Omega_{G\times\mathfrak{g}^*}=\hbox{pr}_1^{*}(\omega_{G\times\mathfrak{g}^*})$, $\Omega_{{\mathfrak g}^{*}}=\hbox{pr}_{{\mathfrak g}^{*}}\omega_{{\mathfrak g}^{*}}$ and $\Omega_U=\hbox{pr}_U^{*}\omega_U$, with $\omega_{G\times\mathfrak{g}^*}$ the left-trivialized symplectic form on $G\times\mathfrak{g}^*$, and $\omega_{{\mathfrak g}^{*}}$ and  $\omega_U$ the canonical symplectic $2$-forms on $T^*\mathfrak{g}^{*}=\mathfrak{g}^{*}\times \mathfrak{g}$ and $T^{*}U$, respectively. Here $\hbox{pr}_1:M\to G\times\mathfrak{g}^{*}$,  $\hbox{pr}_{\mathfrak{g}^{*}}:M\to T^{*}\mathfrak{g}^{*}$ and $\hbox{pr}_{U}:M\to T^{*}U$ are the canonical projections onto the corresponding factors on $M$ and the $*$ notation in the projection maps stands for the pullback of differential forms. Note that locally, these canonical $2$-forms are given by 
\begin{align*}
\Omega_U=&du\wedge dp_u,\\
\Omega_{\mathfrak{g}^{*}}=& d\Pi\wedge dp_{\Pi}\\
(\Omega_{{G\times\mathfrak{g}^*})_{(g, \alpha)}}\left( (\xi_1, \nu_1), (\xi_2, \nu_2)\right)=&-\langle \nu_1, \xi_2\rangle + \langle \nu_2, \xi_1\rangle+\langle\alpha, [\xi_1, \xi_2]\rangle,\label{omega}
\end{align*}
with $(u,p_u)\in T^{*}U$, $(\Pi,p_{\Pi})\in T^{*}{\mathfrak g}^*$, $(g, \alpha)\in G\times {\mathfrak g}^*$,
where $\xi_i\in {\mathfrak g}$ and $\nu_i\in {\mathfrak g}^*$, $i=1,2$.   Note that $\Omega_M$ is presymplectic since $\hbox{Ker }\Omega_M=\displaystyle{\Big{\{}\frac{\partial}{\partial\tau},\frac{\partial}{\partial\dot{u}}\Big{\}}}$.

The equations of motion for this system are given by \begin{equation}\label{intrinsic}i_X\Omega_M=d(\bar{H}_M+\langle p_u,\dot{u}\rangle)\end{equation} for a vector field $X$ on $M$. 

Note that the first constraint submanifold (see Appendix B) is given by \begin{equation}\label{constraint}\frac{\partial\bar{H}_M}{\partial\tau}=0,\quad\frac{\partial\bar{H}_M}{\partial\dot{u}}+p_u=0.\end{equation} Note also that the first condition (i.e., for the external actuators) is the typical condition in optimal control.

In the constraint submanifold which defines equation \eqref{constraint}, solutions of equations \eqref{intrinsic} must satisfy

\begin{align*}
\frac{\partial\bar{H}_M}{\partial\tau}&=0,\quad
  p_u=-\frac{\partial\bar{H}_M}{\partial\dot{u}},\quad 
  \dot{p}_{u}=-\frac{\partial\bar{H}_M}{\partial u},\quad \dot{p}_{\Pi}=-\frac{\partial\bar{H}_M}{\partial \Pi}.\\
\dot{p}_{\xi}&=-T^{*}_{e}\mathcal{L}_g\left(\frac{\partial \bar{H}_M}{\partial g}\right)
+\ad^*_{\xi}p_{\xi},\\ 
\dot{\Pi}&=\bar{F}-I^{-1}(u)\Pi\times\Pi, \quad \dot{g}=gI^{-1}(u)\Pi,
\end{align*} 
where $\xi=g^{-1}\dot{g}$.

That is, \begin{align*}
\frac{\partial C}{\partial \tau}&=p_{\Pi}\frac{\partial \bar{F}}{\partial \tau},\quad
  \frac{d}{dt}\left(\frac{\partial\bar{H}_M}{\partial\dot{u}}\right)-\frac{\partial\bar{H}_M}{\partial u}=0,\\
\dot{p}_{\Pi}&=-\frac{\partial C}{\partial \Pi}+
A_{I(u)}(\Pi)p_{\Pi}-p_{\xi}I^{-1}(u)\\
\dot{p}_{\xi}&=-T^{*}_{e}\mathcal{L}_g\left(\frac{\partial C}{\partial g}\right)
-p_{\Pi}T^{*}_{e}\mathcal{L}_g\left(\frac{\partial \bar{F}}{\partial g}\right)+ad^*_{\xi}p_{\xi}
\\
\dot{\Pi}&=\bar{F}-I^{-1}(u)\Pi\times\Pi, \quad \dot{g}=gI^{-1}(u)\Pi,
\end{align*} where 
$$A_{I(u)}(\Pi)=\begin{pmatrix}
 0&(I_3^{-1}(u)-I_1^{-1}(u))\Pi_3&(I_1^{-1}(u)-I_2^{-1}(u))\Pi_2\\(I_2^{-1}(u)-I_3^{-1}(u))\Pi_3 &0&(I_1^{-1}(u)-I_2^{-1}(u))\Pi_1\\
(I_2^{-1}(u)-I_3^{-1}(u))\Pi_2&(I_3^{-1}(u)-I_1^{-1}(u))\Pi_1&0
\end{pmatrix}$$
and $\xi=g^{-1}\dot{g}$. 
In other words,  
\begin{align}
 \frac{d}{dt}\left(\frac{\partial C}{\partial \dot{u}}\right)&=\frac{\partial C}{\partial u}+p_{\Pi}\left(\frac{\partial }{\partial u}I^{-1}(u)\right)\Pi\times\Pi-p_{\xi}\left(\frac{\partial }{\partial u}I^{-1}(u)\right)\Pi,\\ 
 \frac{\partial C}{\partial \tau}&=p_{\Pi}\frac{\partial \bar{F}}{\partial \tau},\label{eq:tau}\\
\dot{p}_{\Pi}&=-\frac{\partial C}{\partial \Pi}+
A_{I(u)}(\Pi)p_{\Pi}-p_{\xi}I^{-1}(u)\\
\dot{p}_{\xi}&=-T^{*}_{e}\mathcal{L}_g\left(\frac{\partial C}{\partial g}\right)
-p_{\Pi}T^{*}_{e}\mathcal{L}_g\left(\frac{\partial F}{\partial g}\right)+ad^*_{\xi}p_{\xi}\label{pdot}
\\
\dot{\Pi}&=\bar{F}-I^{-1}(u)\Pi\times\Pi, \quad \dot{g}=gI^{-1}(u)\Pi.
\end{align}
Assume that we can isolate the control $\tau$ in Equation (\ref{eq:tau}) as  $\tau=\tau(g,\Pi,u,\dot{u}, p_{\Pi})$ and  $\displaystyle{\det\left(\frac{\partial^{2}C}{\partial\dot{u}\partial\dot{u}}\right)\neq 0}$, then the previous system of differential equations has a unique solution for initial conditions $g(0)$, $\Pi(0)$, $p_{\Pi}(0)$, $p_{\xi}(0)$, $u(0)$, $\dot{u}(0)$.

\subsection{Optimal control problem for foldable drones}\label{sec: foldable}
\begin{figure}[h!]
\begin{center}
\begin{tikzpicture}[line cap=round,line join=round,>=stealth,x=0.5cm,y=0.5cm]
\def\radio{3.5}
\def\brazo{5}
\draw [shift={(0.5*\brazo,0.5*\brazo)},line width=0.25pt,dash pattern=on 1pt off 4pt]  plot[domain=0:6.28,variable=\t]({1*\radio*cos(\t r)+0*\radio*sin(\t r)},{0*\radio*cos(\t r)+1*\radio*sin(\t r)});
\def\radio{0.7}
\draw [shift={(0,0)},line width=2pt]  plot[domain=0:6.28,variable=\t]({1*\radio*cos(\t r)+0*\radio*sin(\t r)},{0*\radio*cos(\t r)+1*\radio*sin(\t r)});
\draw [shift={(\brazo,0)},line width=2pt]  plot[domain=0:6.28,variable=\t]({1*\radio*cos(\t r)+0*\radio*sin(\t r)},{0*\radio*cos(\t r)+1*\radio*sin(\t r)});
\draw [shift={(\brazo,\brazo)},line width=2pt]  plot[domain=0:6.28,variable=\t]({1*\radio*cos(\t r)+0*\radio*sin(\t r)},{0*\radio*cos(\t r)+1*\radio*sin(\t r)});
\draw [shift={(0,\brazo)},line width=2pt]  plot[domain=0:6.28,variable=\t]({1*\radio*cos(\t r)+0*\radio*sin(\t r)},{0*\radio*cos(\t r)+1*\radio*sin(\t r)});
\draw [line width=2pt](0.5*\brazo-0.5,0.5*\brazo-0.5)--(0.5*\brazo-0.5,0.5*\brazo+0.5);
\draw [line width=2pt](0.5*\brazo+0.5,0.5*\brazo-0.5)--(0.5*\brazo+0.5,0.5*\brazo+0.5);
\draw [line width=2pt](0.5*\brazo-0.5,0.5*\brazo-0.5)--(0.5*\brazo+0.5,0.5*\brazo-0.5);
\draw [line width=2pt](0.5*\brazo-0.5,0.5*\brazo+0.5)--(0.5*\brazo,0.5*\brazo+1.0);
\draw [line width=2pt](0.5*\brazo,0.5*\brazo+1.0)--(0.5*\brazo+0.5,0.5*\brazo+0.5);
\draw [line width=2pt](0.5*\brazo-0.5,0.5*\brazo-0.5)--(0,0);
\draw [line width=2pt](0.5*\brazo+0.5,0.5*\brazo+0.5)--(\brazo,\brazo);
\draw [line width=2pt](0.5*\brazo-0.5,0.5*\brazo+0.5)--(0,\brazo);
\draw [line width=2pt](0.5*\brazo+0.5,0.5*\brazo-0.5)--(\brazo,0);
\draw [shift={(0+1.5,0-0.7)},line width=1pt,dash pattern=on 1pt off 4pt]  plot[domain=0:6.28,variable=\t]({1*\radio*cos(\t r)+0*\radio*sin(\t r)},{0*\radio*cos(\t r)+1*\radio*sin(\t r)});
\draw [shift={(\brazo-1.5,0-0.7)},line width=1pt,dash pattern=on 1pt off 4pt]  plot[domain=0:6.28,variable=\t]({1*\radio*cos(\t r)+0*\radio*sin(\t r)},{0*\radio*cos(\t r)+1*\radio*sin(\t r)});
\draw [shift={(\brazo-1.5,\brazo+0.7)},line width=1pt,dash pattern=on 1pt off 4pt]  plot[domain=0:6.28,variable=\t]({1*\radio*cos(\t r)+0*\radio*sin(\t r)},{0*\radio*cos(\t r)+1*\radio*sin(\t r)});
\draw [shift={(0+1.5,\brazo+0.7)},line width=1pt,dash pattern=on 1pt off 4pt]  plot[domain=0:6.28,variable=\t]({1*\radio*cos(\t r)+0*\radio*sin(\t r)},{0*\radio*cos(\t r)+1*\radio*sin(\t r)});
\draw [line width=1pt, dash pattern=on 1pt off 4pt](0.5*\brazo-0.5,0.5*\brazo-0.5)--(0+1.5,0-0.7);
\draw [line width=1pt,dash pattern=on 1pt off 4pt](0.5*\brazo+0.5,0.5*\brazo+0.5)--(\brazo-1.5,\brazo+0.7);
\draw [line width=1pt,dash pattern=on 1pt off 4pt](0.5*\brazo-0.5,0.5*\brazo+0.5)--(0+1.5,\brazo+0.7);
\draw [line width=1pt,dash pattern=on 1pt off 4pt](0.5*\brazo+0.5,0.5*\brazo-0.5)--(\brazo-1.5,0-0.7);
\def\radio{0.7}
\draw [shift={(1.6,1.2)},line width=1pt,->]  plot[domain=0.2:1.8,variable=\t]({-1*\radio*cos(\t r)+0*\radio*sin(\t r)},{0*\radio*cos(\t r)-1*\radio*sin(\t r)});
\draw (0.6,0.8) node[anchor=north west] {$u$};
\draw [shift={(1.6,3.8)},line width=1pt,->]  plot[domain=0.2:1.8,variable=\t]({-1*\radio*cos(\t r)+0*\radio*sin(\t r)},{0*\radio*cos(\t r)+1*\radio*sin(\t r)});
\draw (0.6,5.0) node[anchor=north west] {$u$};
\draw [shift={(3.3,3.8)},line width=1pt,->]  plot[domain=-3.0:-1.4,variable=\t]({-1*\radio*cos(\t r)+0*\radio*sin(\t r)},{0*\radio*cos(\t r)-1*\radio*sin(\t r)});
\draw (3.5,5.0) node[anchor=north west] {$u$};
\draw [shift={(3.3,1.2)},line width=1pt,->]  plot[domain=0.2:1.8,variable=\t]({1*\radio*cos(\t r)+0*\radio*sin(\t r)},{0*\radio*cos(\t r)-1*\radio*sin(\t r)});
\draw (3.5,0.8) node[anchor=north west] {$u$};
\def\radio{0.9}
\draw [shift={(0,0)},line width=1pt,->]  plot[domain=3:5.28,variable=\t]({1*\radio*cos(\t r)+0*\radio*sin(\t r)},{0*\radio*cos(\t r)+1*\radio*sin(\t r)});
\draw [shift={(\brazo,0)},line width=1pt,<-]  plot[domain=4:6.28,variable=\t]({1*\radio*cos(\t r)+0*\radio*sin(\t r)},{0*\radio*cos(\t r)+1*\radio*sin(\t r)});
\draw [shift={(\brazo,\brazo)},line width=1pt,->]  plot[domain=0:2.28,variable=\t]({1*\radio*cos(\t r)+0*\radio*sin(\t r)},{0*\radio*cos(\t r)+1*\radio*sin(\t r)});
\draw [shift={(0,\brazo)},line width=1pt,<-]  plot[domain=1:3.28,variable=\t]({1*\radio*cos(\t r)+0*\radio*sin(\t r)},{0*\radio*cos(\t r)+1*\radio*sin(\t r)});
\draw (\brazo,\brazo+1.7) node[anchor=north west] {$\tau_1$};
\draw (-0.7,\brazo+1.7) node[anchor=north west] {$\tau_2$};
\draw (\brazo,-0.8) node[anchor=north west] {$\tau_4$};
\draw (-0.7,-0.8) node[anchor=north west] {$\tau_3$};
\draw [line width=1pt,->](0.5*\brazo,0.5*\brazo)--(0.5*\brazo,0.5*\brazo+5.5);
\draw [line width=1pt,->](0.5*\brazo,0.5*\brazo)--(0.5*\brazo+5,0.5*\brazo);
\draw [shift={(0.5*\brazo,6.5)},line width=2pt,->]  plot[domain=0.5:2.8,variable=\t]({1*\radio*cos(\t r)+0*\radio*sin(\t r)},{0*\radio*cos(\t r)+1*\radio*sin(\t r)});
\draw [shift={(6,0.5*\brazo)},line width=2pt,->]  plot[domain=-1.14:1.14,variable=\t]({1*\radio*cos(\t r)+0*\radio*sin(\t r)},{0*\radio*cos(\t r)+1*\radio*sin(\t r)});
\def\radio{1.2}
\draw [shift={(0.5*\brazo,0.5*\brazo)},line width=2pt,->]  plot[domain=0:3.14,variable=\t]({1*\radio*cos(\t r)+0*\radio*sin(\t r)},{0*\radio*cos(\t r)+1*\radio*sin(\t r)});
\draw (\brazo+1.7,0.5*\brazo+1.3) node[anchor=north west] {$\phi$};
\draw (0.5*\brazo-1.8,\brazo+3) node[anchor=north west] {$\theta$};
\draw (0.3,0.5*\brazo+1) node[anchor=north west] {$\psi$};
\end{tikzpicture}
\end{center}
\caption{Foldable quadrotor model. In solid line, the vehicle is in ``X'' configuration. In the dashed line, the vehicle is in ``H'' configuration. $u$ denotes the (constant) internal control and $\phi$, $\theta$, $\psi$ denotes roll, pitch and yaw for the quadrotor motion.}\label{drone}
\end{figure}
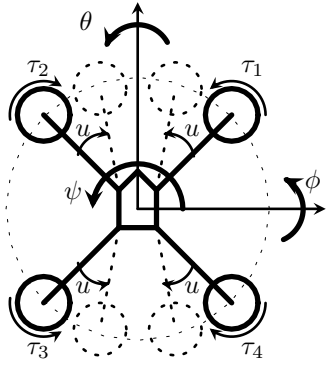

Next, we present the model for the attitude dynamics of a foldable quadrotor considering the impact of the varying nature of the vehicle's inertia. We assume a single rigid body with moment of inertia depending on the (internal controller) angle $u$, assuming constant (see Fig \ref{drone}). 

Consider the Lie group of rotations in the space $SO(3)$, given by $SO(3)=\{R\in GL(3,\mathbb{R})\,|\, R^{T}R=Id, \det(R)=1\}$, where $Id$ denotes the $(3\times 3)$-identity matrix. Denote by $\mathfrak{so}(3)$ its Lie algebra, that is, the set of $3\times 3$ skew-symmetric matrices, given by $\mathfrak{so}(3)=\{\dot{R}_i(0)\,|\,R_i(t)\in SO(3), R_i(0)=Id\}=\{\hat{\omega}_i\in\mathbb{R}^{3\times 3}\,|\,\hat{\omega}_i \hbox{ is skew-symmetric}\}$.

We use the identification of the elements of  the Lie algebra $\mathfrak{so}(3)$ with $\mathbb{R}^{3}$ by $\hat{\omega}_i(t)=\left(
  \begin{array}{ccc}
    0& -\omega^3_i(t) & \omega^2_i(t) \\
    \omega^3_i(t) & 0 & -\omega^1_i(t) \\
    -\omega^2_i(t) & \omega^1_i(t) & 0 \\
  \end{array}
\right)\simeq\omega_i$ with $\omega_i=(\omega_i^{1},\omega_{i}^{2},\omega_i^{3})\in\mathbb{R}^{3}$  and where $\hat{\cdot}:\mathbb{R}^{3}\to\mathfrak{so}(3)$ denotes the isomorphism between $\mathbb{R}^{3}$ and skew-symmetric matrices.

The model is based on the controlled Euler-Poincar\'e equations \cite{Zhai18}, 
\begin{align*}
 I_1(u)\dot{\omega}_1=&(I_2(u)-I_3)\omega_2\omega_3+F_1(u,\tau)-\frac{d}{dt}I_1(u)\omega_1\\
  I_2(u)\dot{\omega}_2=&(I_3-I_1(u))\omega_1\omega_3+F_2(u,\tau)-\frac{d}{dt}I_2(u)\omega_2,\\
 I_3\dot{\omega}_3=&(I_1(u)-I_2(u))\omega_1\omega_2+F_3(u,\tau),
\end{align*} where $(\omega_1,\omega_2,\omega_3)\in\mathbb{R}^{3}$ are the angular velocities of the rigid-body, related to the roll, pitch, and yaw angles $\phi$, $\theta$, $\psi$, respectively. If $\tau=(\tau_1,\tau_2,\tau_3,\tau_4)\in\mathbb{R}^4$ is a vector containing the angular velocities of each rotor (see Fig \ref{drone}), the (external) control signals $F_i$, $i=1,2,3$, are given by 
\begin{align*}
   F_1(u,\tau)=&l\sin(u)\kappa_1(-\tau_1+\tau_2+\tau_3-\tau_4),\\
   F_2(u,\tau)=&l\cos(u)\kappa_1(\tau_1+\tau_2-\tau_3-\tau_4),\\
   F_3(u,\tau)=&l\kappa_2(\tau_1-\tau_2+\tau_3-\tau_4).
\end{align*}With $l$ representing the length of the vehicle's arms, $\kappa_1$ and $\kappa_2$ denoting constants reliant on the motor/propeller characteristics, and where we can analyze the moment of inertia by considering the symmetry of the vehicle. This evaluation accounts for two primary components: the vehicle's center ($I_c$) and its four motors. As such, we derive the following moment of inertia expressions:
\begin{align}I_1(u) &= I_c + 4l^2 \sin^2(u) m_{},\label{i1}\\
I_2(u) &= I_c + 4l^2 \cos^2(u) m_{},\label{i2}\\
I_3 &= I_c + 4l^2 m_{},\label{i3}\end{align} where $m_{}$ is the mass of each motor. 

Using the reduced Legendre transformation $\Pi_i=I_i(u)\omega_i$, the corresponding forced controlled Lie-Poisson equations are given by
\begin{align*}
\dot{\Pi}_1=&\frac{(I_2(u)-I_3)}{I_3I_2(u)}\Pi_2\Pi_3+F_1(u,\tau),\\
 \dot{\Pi}_2=&\frac{(I_3-I_1(u))}{I_3I_1(u)}\Pi_3\Pi_1+F_2(u,\tau),\\
\dot{\Pi}_3=&\frac{(I_1(u)-I_2(u))}{I_1(u)I_2(u)}\Pi_1\Pi_2+F_3(u,\tau).
\end{align*}  

\begin{example}
    For a fixed angle $u$ and a time interval $[t_0,t]$, assuming $\omega_1(z) = \omega_3(z) = 0$, for every $z\in[t_0,t]$, the vehicle's motion will be oriented to the right and/or left. In this specific scenario, the previously stated equations can be further simplified: 
$${\omega}_2(t)=\frac{l\cos(u)\kappa_1}{I_c + 4l^2 \cos^2(u) m_{}}\int_{t_0}^t(\tau_1+\tau_2-\tau_3-\tau_4)_{(z)}dz+\omega_2(t_0).$$
   Considering a given angular velocity trajectory $\omega_2(\bar{t})$ with $\bar{t}\in[t_0,t]$, the optimal angle $u$ can be selected to minimize $\|\tau\|$ and consequently minimize the power consumption. This strategic choice of angle aims to enhance overall efficiency.\hfill$\diamond$
\end{example}

The preceding example highlights how the internal control parameter $u$ significantly influences power consumption. Although the showcased instance assumes a constant angle $u$, real-world applications must also account for the energy necessary to adjust the angle of the vehicle's arms. In such cases, the energy consumption associated with altering the arm angle becomes a factor to consider.

For our optimal control purpose we wish to minimize the cost function \begin{equation}\label{costexample}
C(g, \Pi, u,\dot{u},\tau)=\frac{1}{2}(c_1(\dot{u})^2+c_2||\tau||_{\mathbb{R}^4}^2+c_3||g_d^{-1}g-g^{-1}g_d||_{\mathfrak{so}(3)}^{2}+c_4||\Pi-\Pi_d||^{2}),\end{equation} for a reference state trajectory $g_d,\Pi_d$ and weight constants $c_1,c_2,c_3, c_4> 0$. So, equations \eqref{eq:tau}-\eqref{pdot} are given as follow.

$\bullet$\textit{ Equations for $\tau$:}

\begin{align}
\tau_{1}&=\frac{1}{c_2}\left(-p_{\Pi_{1}}l\kappa_1\sin(u)+p_{\Pi_{2}}\kappa_1l\cos(u)+p_{\Pi_{3}}l\kappa_2\right),\\
\tau_{2}&=\frac{1}{c_2}\left(p_{\Pi_{1}}l\kappa_1\sin(u)+p_{\Pi_{2}}\kappa_1l\cos(u)-p_{\Pi_{3}}l\kappa_2\right),\\
\tau_{3}&=\frac{1}{c_2}\left(p_{\Pi_{1}}l\kappa_1\sin(u)-p_{\Pi_{2}}\kappa_1l\cos(u)+p_{\Pi_{3}}l\kappa_2\right),\\
\tau_{4}&=\frac{-1}{c_2}\left(p_{\Pi_{1}}l\kappa_1\sin(u)+p_{\Pi_{2}}\kappa_1l\cos(u)+p_{\Pi_{3}}l\kappa_2\right).\end{align} Note that we can isolate the external torques as a function of the internal torques, as it is required to stablish regularity of the system. Note that also the regularity condition is fulfilled since $2c_1\neq 0$. 

$\bullet$\textit{ Equations for $p_{\Pi}$:}

\begin{align}
    \dot{p}_{\Pi_1}&=p_{\Pi_2}(I_2^{-1}(u)-I_1^{-1}(u))\Pi_2+p_{\Pi_3}(I_1^{-1}(u)-I_3^{-1})\Pi_3-p_{\omega_1}I_1^{-1}(u)-c_4(\Pi_1-\Pi_1^{d}),\\
     \dot{p}_{\Pi_2}&=p_{\Pi_1}(I_1^{-1}(u)-I_2^{-1}(u))\Pi_1-p_{\Pi_3}(I_2^{-1}(u)-I_3^{-1})\Pi_3-p_{\omega_2}I_2^{-1}(u)-c_4(\Pi_2-\Pi_2^{d}),\\
      \dot{p}_{\Pi_3}&=-p_{\Pi_1}(I_1^{-1}(u)-I_3^{-1})\Pi_1+p_{\Pi_2}(I_2^{-1}(u)-I_3^{-1})\Pi_2-p_{\omega_3}I_3^{-1}-c_4(\Pi_3-\Pi_3^{d}),
\end{align} with $I_1(u)$, $I_2(u)$ and $I_3$ as in equations \eqref{i1}-\eqref{i3}.

$\bullet$\textit{ Equations for $p_{\omega}$:}
\begin{align}
\langle \dot{p}_{\omega}, \eta \rangle&=-2c_3\;\hbox{trace}\left(g_d^{-1}g\,\hat{\eta} \,g_d^{-1}g\right)+ \langle p_{\omega}\times \xi, \eta\rangle, \qquad \forall \eta\in {\mathbb R}^3
.
\end{align}

$\bullet$\textit{ Equations for $u$:}

\begin{align}
    \ddot{u}&=p_{\Pi_1}\left(\frac{\partial}{\partial u}I_2^{-1}(u)\Pi_2\Pi_3-\frac{\partial}{\partial u}I_3^{-1}(u)\Pi_3\Pi_2\right)-p_{\omega_1}\frac{\partial}{\partial u}I_1^{-1}(u)\Pi_1\nonumber,\\
    &+p_{\Pi_2}\left(\frac{\partial}{\partial u}I_3^{-1}\Pi_3\Pi_1-\frac{\partial}{\partial u}I_1^{-1}(u)\Pi_1\Pi_3\right)-p_{\omega_2}\frac{\partial}{\partial u}I_2^{-1}(u)\Pi_2,\\
    &+p_{\Pi_3}\left(\frac{\partial}{\partial u}I_1^{-1}(u)\Pi_1\Pi_2-\frac{\partial}{\partial u}I_2^{-1}(u)\Pi_2\Pi_1\right)-p_{\omega_3}\frac{\partial}{\partial u}I_3^{-1}(u)\Pi_3\nonumber.
\end{align}

$\bullet$\textit{ Equations for $\omega$:}

$$\omega_1=I_1^{-1}(u)\Pi_1,\quad\omega_2=I_2^{-1}(u)\Pi_2,\quad\omega_3=I_3^{-1}\Pi_3.$$

$\bullet$\textit{ Equations for $\Pi$:}

\begin{align}
    \dot{\Pi}_1&=l\sin(u)\kappa_1(-\tau_1+\tau_2+\tau_3-\tau_4)-\frac{\partial}{\partial u}I_2^{-1}(u)\Pi_2\Pi_3-\frac{\partial }{\partial u}I_3^{-1}\Pi_3\Pi_2,\\
    \dot{\Pi}_2&=l\cos(u)\kappa_1(\tau_1+\tau_2-\tau_3-\tau_4)-\frac{\partial}{\partial u}I_3^{-1}\Pi_3\Pi_1+\frac{\partial }{\partial u}I_1^{-1}(u)\Pi_1\Pi_3,\\
    \dot{\Pi}_3&=\kappa_2(\tau_1-\tau_2+\tau_3-\tau_4)-\frac{\partial}{\partial u}I_1^{-1}(u)\Pi_1\Pi_2-\frac{\partial}{\partial u}I_2^{-1}(u)\Pi_2\Pi_1.
\end{align}

One can easily observe that we can not transform through computations the problem to a Hamiltonian function in a symplectic space and construct symplectic integrators. In contrast, in the following section, we will construct variational integrators to integrate numerically the necessary conditions for optimality. 
\section{Variational integrators on Lie groups for control-dependent systems}

In this section, we aim to study the construction of geometric integrators for systems with control-dependent inertia. We first review the discrete Lie-Poisson equations and further develop the numerical integrator for our control-dependent system.

\subsection{Variational integrators on Lie groups}\label{sec: varint} 

Given the set $$\mathcal{T}=\{t_k\in\mathbb{R}^{+},\, t_{k}=kh\mid k=0,\ldots,N\},$$ $Nh=T$, with $T$ fixed, a discrete trajectory for a mechanical system evolving on $G$ is determined by a set of $N+1$ points equally spaced in time, $g_{0:N}=\{g_0,\ldots,g_{N}\}$, where $g_k\simeq g(kh)\in G$, with $h=T/N$ the time step. The path between two adjacent points $g_k$ and $g_{k+1}$ must be given by a curve lying on a manifold. To construct such a curve we make use of a retraction map $\mathcal{R}:\mathfrak{g}\to G$.

A \textit{retraction map} $\mathcal{R}:\mathfrak{g}\to G$ is an 
analytic local diffeomorphism assigning a neighborhood $\mathcal{O}\subset\mathfrak{g}$
of $0\in {\mathfrak g}$ to a neighborhood of the identity $e\in G$.


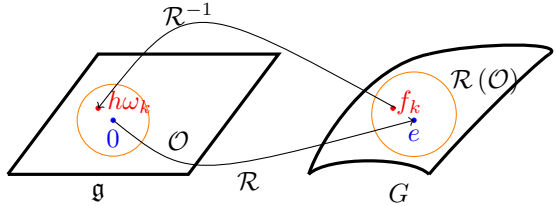
\begin{figure}[htb!] 
\begin{center}
\begin{tikzpicture}[scale=0.8]
\draw[very thick] (4,0) -- (7,0)--(8.5,2)--(5.5,2)--(4,0);
\draw [orange]  (5.75,0.9) circle (17pt);
\filldraw [blue]  (5.75,0.9) circle (1pt) node[below]{$0$};
\filldraw [red] (5.5,1.1) circle(1pt);
\filldraw [red](5.5,1.2)node[right]{$h\omega_k$};
\draw (6.5,0.55) node[right]{$\mathcal{O}$};
\draw (5.5,0) node[below]{$\mathfrak{g}$};
\draw[very thick] (9,0) .. controls (9.2,0.6) and (10,1.8) .. (11,2);
\draw[very thick] (9,0) .. controls (9.3,0.3) and (10.7,0.3) .. (11,0);
\draw[very thick] (11,0) .. controls (11.5,0.6) and (12.7,1.8) .. (13,2);
\draw[very thick] (11,2) .. controls (11.8,2.2) and (13.3,2.2) .. (13,2);
\draw [orange]  (10.75,1) circle (20pt);
\filldraw [blue]  (10.75,0.9) circle (1pt) node[below]{$e$};
\filldraw [red] (10.4,1.1) circle(1pt); 
\filldraw [red](10.3,1.2)node[right]{$f_{k}$};
\draw (11.2,1.55) node[right]{$\mathcal{R}\left(\mathcal{O}\right)$};
\draw (10.5,0) node[below]{$G$};
\draw [->]  (5.75,0.9) .. controls (7,-0.1) and (7,-0.1) .. (10.75,0.9);
\draw [<-](5.5,1.1) .. controls (7,3) and (7,3) .. (10.45,1.1);
\draw (7,3) node [below]{$\mathcal{R}^{-1}$};
\draw (8,0.2) node [below]{$\mathcal{R}$};
\end{tikzpicture}\caption{Retraction map.}\label{retractionfigure}
\end{center}
\end{figure}

The retraction map  (see Figure \ref{retractionfigure}) is used to express small discrete changes in the
group configuration through unique Lie algebra elements given by $$\omega_k=\mathcal{R}^{-1}(g_{k}^{-1}g_{k+1})/h,$$ where
$\omega_k\in\mathfrak{g}$ (see \cite{Rabee} and \cite{KM} for
further details). That is, if $\omega_{k}$ were regarded as an
average velocity between $g_{k}$ and $g_{k+1}$, then $\mathcal{R}$ is an
approximation to the integral flow of the dynamics. The difference
$f_{k}:=g_{k}^{-1}g_{k+1}\in G$, which is an element of a nonlinear
space, can now be represented by the vector $\omega_{k},$ in order to
enable unconstrained optimization in the linear space $\mathfrak{g}$ for optimal control purposes.

For the derivation of the
discrete equations of motion, the {\it right trivialized}
tangent retraction map will be used. It is the  function
$d\mathcal{R}:\mathfrak{g}\times\mathfrak{g}\rightarrow\mathfrak{g}$ given by \begin{equation}\label{righttrivialization}
T\mathcal{R}(\xi)\cdot\eta=TR_{\mathcal{R}(\xi)}d\mathcal{R}_{\xi}(\eta),
\end{equation}
where $\eta\in\mathfrak{g}$ and $R:G\to G$ the right translation on $G$ (see \cite{Rabee} for the derivation of such a map). Here we use the following notation,
$d\mathcal{R}_\xi:=d\mathcal{R}(\xi):\mathfrak{g}\rightarrow\mathfrak{g}.$ The
function $d\mathcal{R}$ is linear, but only on one argument. 

The natural choice of a retraction map is the exponential map at
the identity $e$ of the group $G,$ $\hbox{exp}_{e}:\mathfrak{g}\rightarrow
G$. Recall that, for a finite-dimensional Lie group, $\hbox{exp}_{e}$ is
locally a diffeomorphism and gives rise to a natural chart
\cite{Hair}. Then, there exists a neighborhood $U$ of $e\in G$ such
that $\hbox{exp}_{e}^{-1}:U\rightarrow \hbox{exp}_{e}^{-1}(U)$ is a local
$\mathcal{C}^{\infty}-$diffeomorphism. A chart at $g\in G$ is given
by $\Psi_{g}=\hbox{exp}_{e}^{-1}\circ L_{g^{-1}}$. In general, it is not easy to work with the exponential map since the differential of the exponential map involves power series expansions with iterated Lie-brackets. In
consequence it will be useful to use a different retraction map.
More concretely, the Cayley map, which is usually used in numerical integration with matrix Lie-groups configurations (see \cite{Rabee, Hair, KM} for further
details), will provide to us a proper framework in the application shown
in the paper.

Next, we consider the discrete Lagrangian $\ell_{d}\colon G\to\mathbb{R}$  as an approximation of $\ell\colon  \mathfrak{g}\to\mathbb{R}$ to construct variational integrators in
the same way as in discrete mechanics \cite{mawest}. For a given $h>0$ we define
$\ell_d:G\to\mathbb{R}$ as an approximation of the action integral for $\ell$
along each discrete segment between $g_k$ and $g_{k+1}$, that is, $$\int_{kh}^{(k+1)h}\ell(\omega)\,dt\simeq h\ell_{d}\left(\kappa(g_k,g_{k+1})\right),$$
where $\kappa$ is a function of
$(g_k,g_{k+1})\in G\times G$ which approximate  $\omega(t)$. The discrete Lie-Poisson equations are given by (see \cite{Rabee} for instance) \begin{align}
g_{k+1} =&g_k\mathcal{R}(h\omega_k),\label{reconstructiondiscrete}\\
\left(d\mathcal{R}^{-1}_{(h\omega_k)}\right)^{*}\Pi_k=&\left(d\mathcal{R}^{-1}_{(-h\omega_{k-1})}\right)^{*}\Pi_{k-1}, \label{eqlpdiscrete}\\
\Pi_k=&\frac{\partial \ell_d}{\partial \omega_k},\label{momentumdiscrete}
\end{align} for $k=1,\ldots,N-1$.

\subsection{Variational integrators for the optimal control of controlled inertia rigid bodies}\label{sec4.2}

Next, we will construct a variational integrator for the necessary conditions for optimality in the optimal control problem of controlled inertia rigid bodies given in Section \ref{ocp}.

Consider the discretization of the cost function \eqref{costexample} $C:G\times\mathfrak{g}^{*}\times TU\times\bar{U}\to\mathbb{R}$ as $C_d:G\times G\times\mathfrak{g}^{*}\times\mathfrak{g}^{*}\times U\times U\times\bar{U}\times\bar{U}\to\mathbb{R}$ given by the trapezoidal discretization of $C$, that is, 
\begin{align}
C_d(g_k,g_{k+1},\Pi_k,\Pi_{k+1},u_k,u_{k+1},\tau_k,\tau_{k+1})&=\frac{1}{2}C\left(g_k,\Pi_k,u_k,\frac{u_{k+1}-u_k}{h}, \tau_k\right)\nonumber\\
&+\frac{1}{2}C\left(g_{k+1},\Pi_{k+1},u_{k+1},\frac{u_{k+1}-u_k}{h},\tau_{k+1}\right),\label{discretecost}
\end{align} provided that the set of controls $U$ has a vector space structure. 

The external force $F:G\times\mathfrak{g}^*\times U\times\bar{U}\to \mathfrak{g}^{*}$ is also discretized by using the trapezoidal discretization, that is, $F^k_d:G\times G\times\mathfrak{g}^{*}\times\mathfrak{g}^{*}\times U\times U\times\bar{U}\times\bar{U}\to {\mathfrak{g}^{*}}$ is given by 
\begin{align}
F^k_d(g_k,g_{k+1},\Pi_k,\Pi_{k+1},u_k,u_{k+1},\tau_k,\tau_{k+1})&=\frac{1}{2}F\left(g_k,\Pi_k,u_k, \tau_k\right)\nonumber\\
&+\frac{1}{2}F\left(g_{k+1},\Pi_{k+1},u_{k+1},\tau_{k+1}\right).\label{discreteforce}
\end{align}
Finally, we consider the trapezoidal discretization for the equations \eqref{kinematics1}-\eqref{kinematics2}, that is, 
\begin{align}
0=&\frac{\Pi_{k+1}-\Pi_k}{h}-F_d^ k(c_{k,k+1})+\frac{1}{2}I^{-1}(u_k)\Pi_k\times\Pi_k+\frac{1}{2}I^{-1}(u_{k+1})\Pi_{k+1}\times\Pi_{k+1},\label{kinematics1disc}\\
0=&\mathcal{R}^{-1}(g_{k}^{-1}g_{k+1})-\frac{h}{2}I^{-1}(u_k)\Pi_k-\frac{h}{2}I^{-1}(u_{k+1})\Pi_{k+1},\label{kinematics2disc}
\end{align} where $c_{k,k+1}:=(g_k,g_{k+1},\Pi_k,\Pi_{k+1},u_k,u_{k+1},\tau_k,\tau_{k+1})$ and $\mathcal{R}:\mathfrak{g}\to G$ is a retraction map.

To derive the discrete equations of motion we first denote by 

\begin{align*}
&\Phi_{h}^{k}(g_k,g_{k+1},u_{k},u_{k+1},\Pi_k,\Pi_{k+1}):=\mathcal{R}^{-1}(g_{k}^{-1}g_{k+1})-\frac{h}{2}I^{-1}(u_k)\Pi_k-\frac{h}{2}I^{-1}(u_{k+1})\Pi_{k+1}\\
&\widetilde{F}^k_d(c_{k,k+1}):=\frac{1}{2}I^{-1}(u_k)\Pi_k\times \Pi_k+\frac{1}{2}I^{-1}(u_{k+1})\Pi_{k+1}\times \Pi_{k+1}-F_d^k(c_{k,k+1}),
\end{align*} and we consider the extended discrete cost function \begin{align}
\widetilde{C}_d(c_{k,k+1},\lambda_k,\mu_k)=&\lambda_k\left(\frac{\Pi_{k+1}-\Pi_k}{h}+\widetilde{F}_k(c_{k,k+1})\right)\\&+\langle\mu_k,\Phi_{h}^{k}(u_k,u_{k+1},\Pi_k,\Pi_{k+1})\rangle -C_d(c_{k,k+1}).\nonumber \end{align} where $\mu_k\in \mathfrak{g}^{*}$ and $\lambda_k\in\mathfrak{g}$ are the Lagrange multipliers.

Consider the discrete action $$S_d=\sum_{k=0}^{N-1}\widetilde{C}_d(g_k,g_{k+1},\Pi_k,\Pi_{k+1},u_k,u_{k+1},\tau_k,\tau_{k+1},\lambda_k,\mu_k).$$ Taking variations with fixed endpoints and   shifting the order of the indexes as usual in discrete variational calculus, we obtain the following set of algebraic difference equations
\begin{align}
0&=-D_{1,\tau_k^i}C_d^{k}-D_{2,\tau_k^i}C_d^{k-1}+\langle\lambda_{j,k},D_{1,\tau_k^i}F_{d}^{j,k}\rangle+\langle\lambda_{j,k-1},D_{2,\tau_k^i}F_d^{j,k-1}\rangle,\label{eqq1}\\
0&=-D_{1,\Pi_k^i}C_d^k-D_{2,\Pi_k^i}C_d^{k-1}-\frac{\lambda_k-\lambda_{k-1}}{h}+\lambda_kD_{1,\Pi_k^i}\widetilde{F}_d^k+\lambda_{k-1}D_{2,\Pi_k^i}\widetilde{F}_d^{k-1}\nonumber\\
&-\mu_kD_{1,\Pi_k^i}\Phi_{h}^{k}-\mu_{k-1}D_{2,\Pi_k^i}\Phi_{h}^{k-1},\label{eqq2}\\
0&=-T^{*}_{e}\mathcal{L}_{g_k}(D_{1,g_k}C_d^k)-T^{*}_{e}\mathcal{L}_{g_k}(D_{2,g_k}C_d^{k-1})-\lambda_k T^{*}_{e}\mathcal{L}_{g_k}F_d^k\nonumber\\
&-\lambda_{k-1}T^{*}_{e}\mathcal{L}_{g_k}F_d^{k-1}-T^{*}_eR_{f_k}\mu_k+T^{*}_{e}\mathcal{L}_{f_{k-1}}\mu_{k-1},\label{eqq3}\\
0&=-D_{1,u_k}C_d^k-D_{2,u_k}C_d^{k-1}+\langle\lambda_k,D_{1,u_k}\widetilde{F}_{d}^{k}\rangle+\langle\lambda_{k-1},D_{2,u_k}\widetilde{F}_d^{k-1}\rangle\nonumber\\&-\langle\mu_k,D_{1,u_k}\Phi_{h}^{k}\rangle-\langle\mu_{k-1},D_{2,u_k}\Phi_{h}^{k-1}\rangle\label{eqq4}.
\end{align} together with equations \eqref{kinematics1disc}-\eqref{kinematics2disc} where $f_k=g_k^{-1}g_{k+1}$. Here we have employed the following notation: for a function $P$, $P^k$ denotes $P^{k}:=P(c_{k,k+1})$. $D_{j,s}P$ denotes the partial derivative of $P$ in the $j$-position of $P$ where appears the variable $s$. In addition $\lambda_{j,k}$ and $F^{j,k}_d$ means the $j$-component of $\lambda_k$ and $F_d^k$ respectively, where $j=1,\ldots,|F^k_d|$, where we understand $F_d^{k}$ as  vector with components $F_d^{j,k}$.

Assume we can isolate $\tau_{k+1}$ in terms of $(g_k,g_{k+1},\Pi_k,\Pi_{k+1},u_k,u_{k+1},\tau_k,\lambda_k,\mu_k)$ in equation \eqref{eqq1}, if in addition the matrix $\left(-D_{12,u_k}C_d^k+\lambda_kD_{12,u_k}\widetilde{F_d^k}-\mu_kD_{12}\Phi_{h}^{k}\right)$ is non-singular (i.e., invertible), by a direct application of the implicit function theorem, we deduce that, in a neighborhood $\mathcal{U}_d$ of the point\newline $c=(g,\bar{g},\Pi,\bar{\Pi},u,\bar{u},\tau,\bar{\tau},\lambda,\mu)\in 2(G\times \mathfrak{g}^{*}\times U\times \bar{U})\times\mathfrak{g}\times G$, exists a unique local application \[
\begin{array}{rrcl}
\Upsilon_d:& {\mathcal U}_d\longrightarrow& {\mathcal U}_d\\
       & c \longmapsto&\bar{c}\; ,
\end{array} 
\]with $\bar{c}=(\bar{g},\hat{g},\bar{\Pi},\hat{\Pi},\bar{u},\hat{u},\bar{\tau},\hat{\tau},\bar{\lambda},\bar{\mu})\in 2G\times 2\mathfrak{g}^{*}\times 2U\times 2\bar{U}\times\mathfrak{g}\times G$, such that for all solutions $(c_{0,1},c_{1,2},\ldots,c_{N-1,N},\lambda_0,\ldots,\lambda_N,\mu_0,\ldots,\mu_N)$ of the discrete equations \eqref{eqq1}-\eqref{eqq4} together with \eqref{kinematics1disc} and \eqref{kinematics2disc}, we have that locally $$\Upsilon_d(c_{k,k+1},\lambda_k,\mu_k)=(c_{k+1,k+2},\lambda_{k+1},\mu_{k+1}).$$

Note also that the condition $\left(-D_{12,u_k}C_d^k+\lambda_kD_{12,u_k}\widetilde{F_d^k}-\mu_kD_{12}\Phi_{h}^{k}\right)$ 
is related with the continuous-time regularity condition by 
$$\displaystyle{\left(-D_{12,u_k}C_d^k+\lambda_kD_{12,u_k}\widetilde{F_d^k}-\mu_kD_{12}\Phi_{h}^{k}\right)\simeq -\frac{1}{h}\left(\frac{\partial^{2}C}{\partial\dot{u}\partial\dot{u}}\right)+\mathcal{O}(1)}\; .$$ 

\section{Variational integrators on Lie groups for the optimal control of  foldable drones}\label{sec5}

Next we describe the optimal control problem to reach a desired pose for a foldable UAV.

We consider the following discrete-time optimal control problem 

$$\min_{\{g_k,u_k,\Pi_k,\tau_k\}}\sum_{k=0}^{N-1}C_d(g_k,g_{k+1},\Pi_k,\Pi_{k+1}, u_k,u_{k+1},\tau_k,\tau_{k+1})$$ with $g(0)=g_0$, $\Pi(0)=\Pi_0$, $u(0)=u_0$, $g(T)=g_N$, $\Pi(T)=\Pi_N$, $u(T)=u_N$ fixed boundary conditions, and subject to 
\begin{align}&\Phi_1(\Pi_k,\Pi_{k+1},u_k,u_{k+1},\tau_k,\tau_{k+1})=\frac{\Pi_{k+1}^{1}-\Pi_k^1}{h}-\frac{I_2(u_k)-I_3}{2I_3I_2(u_k)}\Pi_{k}^{2}\Pi_k^{3}\label{cons1}\\&\qquad\qquad\qquad\qquad-\frac{I_2(u_{k+1})-I_3}{2I_3I_2(u_{k+1})}\Pi_{k+1}^{2}\Pi_{k+1}^{3}-\frac{1}{2}F_1(u_k,\tau_k)-\frac{1}{2}F_1(u_{k+1},\tau_{k+1})=0, \nonumber\\ 
&\Phi_2(\Pi_k,\Pi_{k+1},u_k,u_{k+1},\tau_k,\tau_{k+1})=\frac{\Pi_{k+1}^{2}-\Pi_k^2}{h}-\frac{I_3-I_1(u_k)}{2I_3I_1(u_k)}\Pi_{k}^{3}\Pi_k^{1}\label{cons2}\\&\qquad\qquad\qquad-\frac{I_3-I_1(u_{k+1})}{2I_3I_1(u_{k+1})}\Pi_{k+1}^{3}\Pi_{k+1}^{1}-\frac{1}{2}F_2(u_k,\tau_k)-\frac{1}{2}F_2(u_{k+1},\tau_{k+1})=0,\nonumber \\
&\Phi_3(\Pi_k,\Pi_{k+1},u_k,u_{k+1},\tau_k,\tau_{k+1})=\frac{\Pi_{k+1}^{3}-\Pi_k^3}{h}-\frac{I_1(u_k)-I_2(u_k)}{2I_1(u_k)I_2(u_k)}\Pi_{k}^{1}\Pi_k^{2}\label{cons3}\\&\qquad\qquad-\frac{I_1(u_{k+1})-I_2(u_{k+1})}{2I_1(u_{k+1})I_2(u_{k+1})}\Pi_{k+1}^{1}\Pi_{k+1}^{2}-\frac{1}{2}F_3(u_k,\tau_k)-\frac{1}{2}F_3(u_{k+1},\tau_{k+1})=0,\nonumber \\
&\Phi_4(\Pi_k,\Pi_{k+1},u_k,u_{k+1})=
\mathcal{R}^{-1}(g_k^{-1}g_{k+1})-\frac{h}{2}I^{-1}(u_k)\Pi_k-\frac{h}{2}I^{-1}(u_{k+1})\Pi_{k+1}=0,\label{cons4}
\end{align} where $$I^{-1}(u_k)\Pi_k:=\left(\begin{array}{ccc}
    I_1^{-1}(u_k)& 0 & 0 \\
   0 & I_2^{-1}(u_k) & 0 \\
   0 & 0 & I_3^{-1}(u_k) \\
  \end{array}\right)\left(\begin{array}{c}
    \Pi_k^{1}\\
   \Pi_k^{2} \\
\Pi_k^{3} \\
  \end{array}\right),$$ and where the discrete cost function is given by 
\begin{align*}
   &C_d(g_k,g_{k+1},\Pi_k,\Pi_{k+1},u_k,u_{k+1},\tau_k,\tau_{k+1})=\frac{c_1}{2h}(u_{k+1}-u_k)^2+ \frac{c_2}{4}\sum_{i=1}^{4}(\tau_k^{i})^2+\frac{c_2}{4}\sum_{i=1}^{4}(\tau_{k+1}^{i})^2\\&+\frac{c_3}{4}\hbox{Tr}\left((g_d^{-1}(kh)g_k-g_k^{-1}g_d(kh))^{T}(g_d^{-1}(kh)g_k-g_k^{-1}g_d(kh))\right)\\&+\frac{c_3}{4}\hbox{Tr}\left((g_d^{-1}((k+1)h)g_{k+1}-g_{k+1}^{-1}g_d((k+1)h))^{T}(g_d^{-1}((k+1)h)g_{k+1}-g_{k+1}^{-1}g_d((k+1)h))\right)\\
   &+\frac{c_4}{4}||\Pi_k-\Pi_d(kh)||^{2}+\frac{c_4}{4}||\Pi_{k+1}-\Pi_d((k+1)h)||^{2})
\end{align*}

We consider the extended cost function $\widetilde{C}_d$ by including the constraints $\Phi_i$, $i=1,\ldots,4$ together with some Lagrange multipliers $\lambda_k^{1}$, $\lambda_k^{2}$, $\lambda_k^{3}$, $\mu_k$ into the cost functional, that is, $$\widetilde{C}_d=\langle\lambda_k^{1},\Phi_1\rangle+\langle\lambda_k^{2},\Phi_2\rangle+\langle\lambda_k^{3},\Phi_3\rangle+\langle\mu_k,\Phi_4\rangle-C_d$$ with $\lambda_k^{i}\in\mathbb{R}$ for $i=1,2,3$ and $\mu_k\in G$. So, equations \eqref{eqq1}-\eqref{eqq4}, together with the constraints \eqref{cons1}-\eqref{cons4}, are given as follow:

$\bullet$\textit{ Equations for $\tau_k$:}
\begin{align}
\tau_{k}^{1}&=-\frac{1}{c_2}\left(-\lambda_{1,k}l\kappa_1\sin(u_k)+\lambda_{2,k}\kappa_1l\cos(u_k)+\lambda_{3,k}l\kappa_2\right),\\
\tau_{k}^{2}&=-\frac{1}{c_2}\left(\lambda_{1,k}l\kappa_1\sin(u_k)+\lambda_{2,k}\kappa_1l\cos(u_k)-\lambda_{3,k}l\kappa_2\right),\\
\tau_{k}^{3}&=-\frac{1}{c_2}\left(\lambda_{1,k}l\kappa_1l\sin(u_k)-\lambda_{2,k}\kappa_1l\cos(u_k)+\lambda_{3,k}l\kappa_2\right),\\
\tau_{k}^{4}&=-\frac{1}{c_2}\left(\lambda_{1,k}l\kappa_1l\sin(u_k)+\lambda_{2,k}\kappa_1l\cos(u_k)+\lambda_{3,k}l\kappa_2\right).\end{align} Note that we can isolate the external torques as a function of the internal torques, as it is required in the continuous-time situation to guarantee regularity of the system.

$\bullet$\textit{ Equations for $\lambda_k$:} For the derivation of the equations for $\Pi_k$ we employ the cayley map on $SO(3)$ as retraction map \cite{muller2021review}. The Cayley transform is the map $\hbox{cay}:\hat{x}\in\mathfrak{so}(3)\mapsto (I-\hat{x})^{-1}(I+\hat{x})\in SO(3).$
Indeed $\hbox{cay}(\hat{x})^{T}\hbox{cay}(\hat{x})=I$. We then have the formulas 

\begin{align}
    \hbox{cay}(\hat{x})&=I+2\sigma\hat{x}+2\sigma\hat{x}^{2}\\
    \hbox{dcay}(\hat{x})&=2\sigma(I+\hat{x}),
\end{align}where $\sigma:=\frac{1}{1+||x||^2}$. Then the equations for $\Pi_k$ have the following expressions

\begin{align}
\frac{\lambda_k^1-\lambda_{k-1}^1}{h}=&-\lambda_k^{2}\left(\frac{I_3+I_1(u_k)}{2I_3I_1(u_k)}\right)\Pi_k^3-\lambda_{k}^{3}\left(\frac{I_2(u_k)+I_1(u_k)}{2I_2(u_k)I_1(u_k)}\right)\Pi_k^2\nonumber\\&-\lambda_{k-1}^{2}\left(\frac{I_3+I_1(u_k)}{2I_3I_1(u_k)}\right)\Pi_k^3-\lambda_{k-1}^{3}\left(\frac{I_2(u_k)+I_1(u_k)}{2I_2(u_k)I_1(u_k)}\right)\Pi_k^2\\
&-\mu_k\frac{2}{1+||\tilde{x}||}(I+\hat{\tilde{x}})-\mu_{k-1}\frac{2}{1+||\tilde{y}||}(I+\hat{\tilde{y}})\\
&-c_4(\Pi^1_k-\Pi^1_d(kh)),
,\nonumber\\
\frac{\lambda_k^2-\lambda_{k-1}^2}{h}=&\lambda_k^{2}\left(\frac{I_3+I_2(u_k)}{2I_3I_2(u_k)}\right)\Pi_k^3-\lambda_{k}^{3}\left(\frac{I_2(u_k)+I_1(u_k)}{2I_2(u_k)I_1(u_k)}\right)\Pi_k^1\nonumber\\&-\lambda_{k-1}^{2}\left(\frac{I_3+I_2(u_k)}{2I_3I_2(u_k)}\right)\Pi_k^3-\lambda_{k-1}^{3}\left(\frac{I_2(u_k)+I_1(u_k)}{2I_2(u_k)I_1(u_k)}\right)\Pi_k^1\\
&-\mu_k\frac{2}{1+||\tilde{x}||}(I+\hat{\tilde{x}})-\mu_{k-1}\frac{2}{1+||\tilde{y}||}(I+\hat{\tilde{y}})\\
&-c_4(\Pi^2_k-\Pi^2_d(kh)),\nonumber\\
\frac{\lambda_k^3-\lambda_{k-1}^3}{h}=&\lambda_k^{2}\left(\frac{I_3+I_2(u_k)}{2I_3I_2(u_k)}\right)\Pi_k^2-\lambda_{k}^{3}\left(\frac{I_3+I_1(u_k)}{2I_3I_1(u_k)}\right)\Pi_k^1\nonumber\\&-\lambda_{k-1}^{2}\left(\frac{I_3+I_2(u_k)}{2I_3I_2(u_k)}\right)\Pi_k^2-\lambda_{k-1}^{3}\left(\frac{I_3+I_1(u_k)}{2I_3I_1(u_k)}\right)\Pi_k^1\\
&-\mu_k\frac{2}{1+||\tilde{x}||}(I+\hat{\tilde{x}})-\mu_{k-1}\frac{2}{1+||\tilde{y}||}(I+\hat{\tilde{y}})\\
&-c_4(\Pi^3_k-\Pi^3_d(kh)),\nonumber
\end{align} where 

\begin{align}
\tilde{x}=&\frac{h}{2}I^{-1}(u_k)\Pi_k+\frac{h}{2}I^{-1}(u_{k+1})\Pi_{k+1}\\ \tilde{y}=&\frac{h}{2}I^{-1}(u_{k-1})\Pi_{k-1}+\frac{h}{2}I^{-1}(u_{k})\Pi_{k}
\end{align} and the norm is a norm on $\mathfrak{so}^{*}(3)$.

$\bullet$\textit{ Equations for $u_k$:}
\begin{align}\label{equ}
    0=&\frac{c_1}{h}(u_{k+1}-2u_k+u_{k-1})-\frac{4\mu_k^{1}hl^2\cos(u_k)m}{(I_c+4l^2\sin^{2}(u_k)m)^2}\Pi_k^{1}-\frac{4\mu_k^{2}hl^2\sin(u_k)m}{(I_c+4l^2\cos^{2}(u_k)m)^2}\Pi_k^{2}\\
    &-\lambda_k^{1}\frac{4\ell^2\cos(u_k)m}{(I_c+4\ell^2\sin^{2}(u_k)m)^2}\Pi_k^{1}\Pi_k^{3}-\lambda_k^{2}\frac{4\ell^2\sin(u_k)m}{(I_c+4\ell^2\cos^{2}(u_k)m)^2}\Pi_k^{2}\Pi_k^{3}\nonumber\\
    &-\lambda_k^{3}\frac{4\ell^2\sin(u_k)m}{(I_c+4\ell^2\cos^{2}(u_k)m)^2}\Pi_k^{1}\Pi_k^{2}-\lambda_{k-1}^{1}\frac{4\ell^2\cos(u_k)m}{(I_c+4\ell^2\sin^{2}(u_k)m)^2}\Pi_k^{1}\Pi_k^{3}-\nonumber\\
    &-\lambda_{k-1}^{2}\frac{4\ell^2\sin(u_k)m}{(I_c+4\ell^2\cos^{2}(u_k)m)^2}\Pi_k^{2}\Pi_k^{3}-\lambda_{k-1}^{3}\frac{4\ell^2\sin(u_k)m}{(I_c+4\ell^2\cos^{2}(u_k)m)^2}\Pi_k^{1}\Pi_k^{2}\nonumber\\
    &-\frac{4\mu_{k-1}^{1}hl^2\cos(u_k)m}{(I_c+4l^2\sin^{2}(u_k)m)^2}\Pi_k^{1}-\frac{4\mu_{k-1}^{2}hl^2\sin(u_k)m}{(I_c+4l^2\cos^{2}(u_k)m)^2}\Pi_k^{2}.\nonumber
\end{align}

$\bullet$\textit{ Equations for $g_k$:}
\begin{align}
0=&2c_3\hbox{Tr}\left(g_d(kh)^{-1}g_k \eta g_d(kh)^{-1}g_k \right)+ \langle \mu_k,T_{f_k}\hbox{cay}^{-1}(T_e R_{f_k}(\eta))\rangle\nonumber\\
&-  \langle \mu_{k-1},T_{f_{k-1}}\hbox{cay}^{-1}(T_e L_{f_{k-1}}(\eta))\rangle ,
\end{align} where $f_k=g_k^{-1}g_{k+1}$ and $\eta\in\mathfrak{so}(3)^{*}$ arbitrary.

\subsection{Simulation results}
This section presents the results of a simulation experiment conducted to evaluate the proposed orientation control strategy on a simulated quadrotor. The specifications of these quadrotor resembles the DJI-F450, a commonly used quadrotor for research purposes, with a weight of approximately 1.6 kg (including a standard battery) and approximate moments of inertia given by $I_x = 0.034$ kg·m\(^2\), $I_y = 0.034$ kg·m\(^2\), and $I_z = 0.056$ kg·m\(^2\). These values were obtained from manufacturer data and typical configurations of the vehicle in nominal arm angle position.

The simulation considered only the rotational dynamics of the quadrotor, neglecting translational motion. In order to illustrate the behavior of the quadrotor, Euler angles are used. These angles represent the orientation of the vehicle with respect to a fixed reference frame and are defined as follows: \textit{yaw} (rotation around the vertical axis, describing the heading of the vehicle), \textit{pitch} (rotation around the lateral axis, describing the nose-up or nose-down attitude), and \textit{roll} (rotation around the longitudinal axis, describing the tilting of the vehicle to the sides).

\begin{figure}[h!]
    \centering
    \includegraphics[width=0.6\textwidth, angle=270]{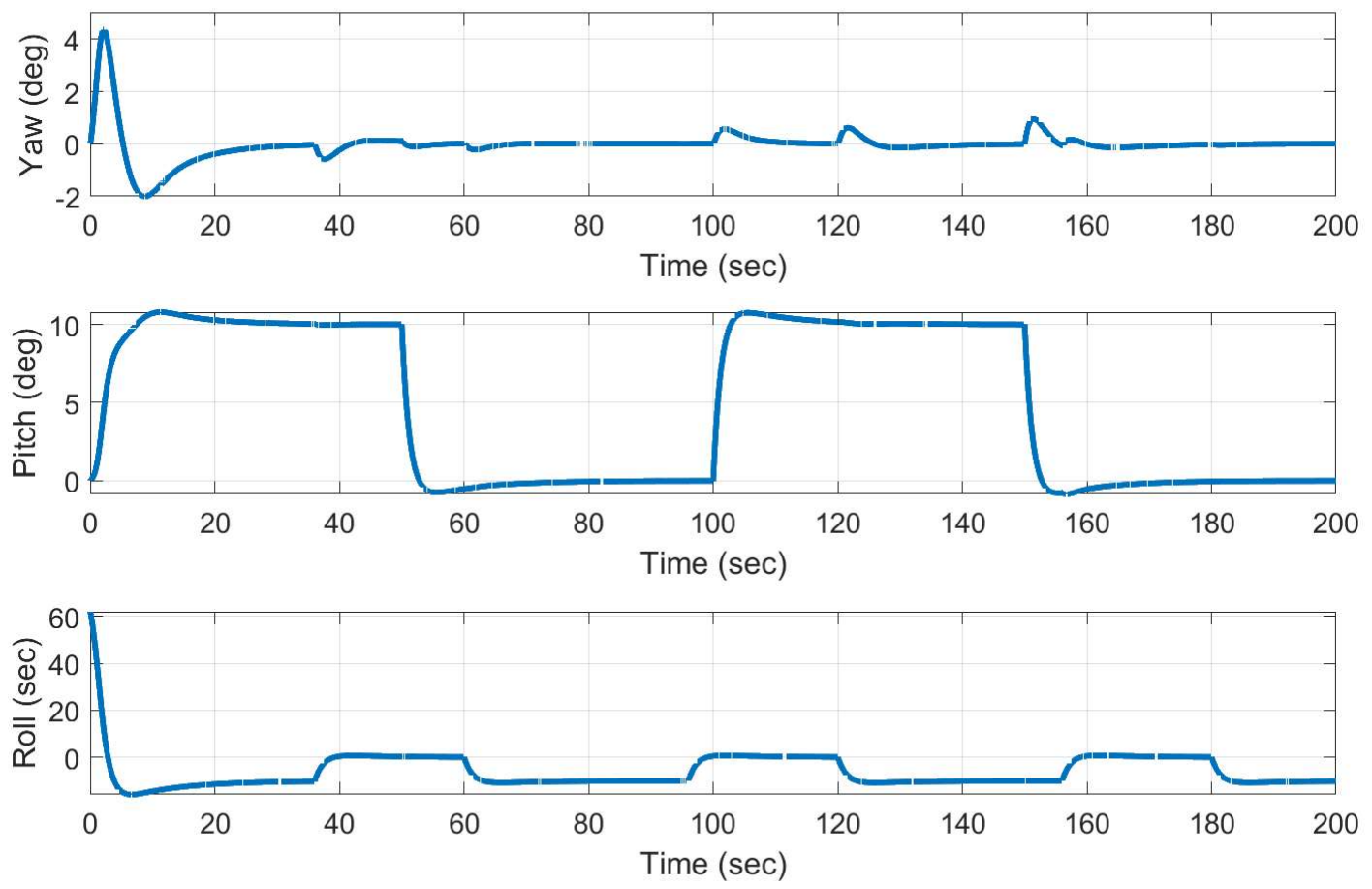}
    \caption{Vehicle attitude (roll, pitch, and yaw) during stabilization and tracking phases.}
    \label{fig:attitude}
\end{figure}

The experiment began with the quadrotor initialized in a highly deviated orientation far from the desired zero angles. Despite these challenging initial conditions, the controller successfully stabilized the vehicle, bringing it to a steady-state orientation in a short amount of time. Figure~\ref{fig:attitude} illustrates this stabilization phase, showing the convergence of the vehicle's attitude (roll, pitch, and yaw) to the reference orientation. The initial values are $p_\Pi(0)=( 0.0752  ,0.0091,0.0049)$, $p_\xi(0)=(0.0227,   0.0027,    0.0020)$, with an initial roll tilt of $1.0821rad$. The weight for the functional are $(c_1,c_2,c_2,c_3)=(0.01,1,1,0.1)$.

\begin{figure}[h!]
    \centering
    \includegraphics[width=0.6\textwidth, angle=270]{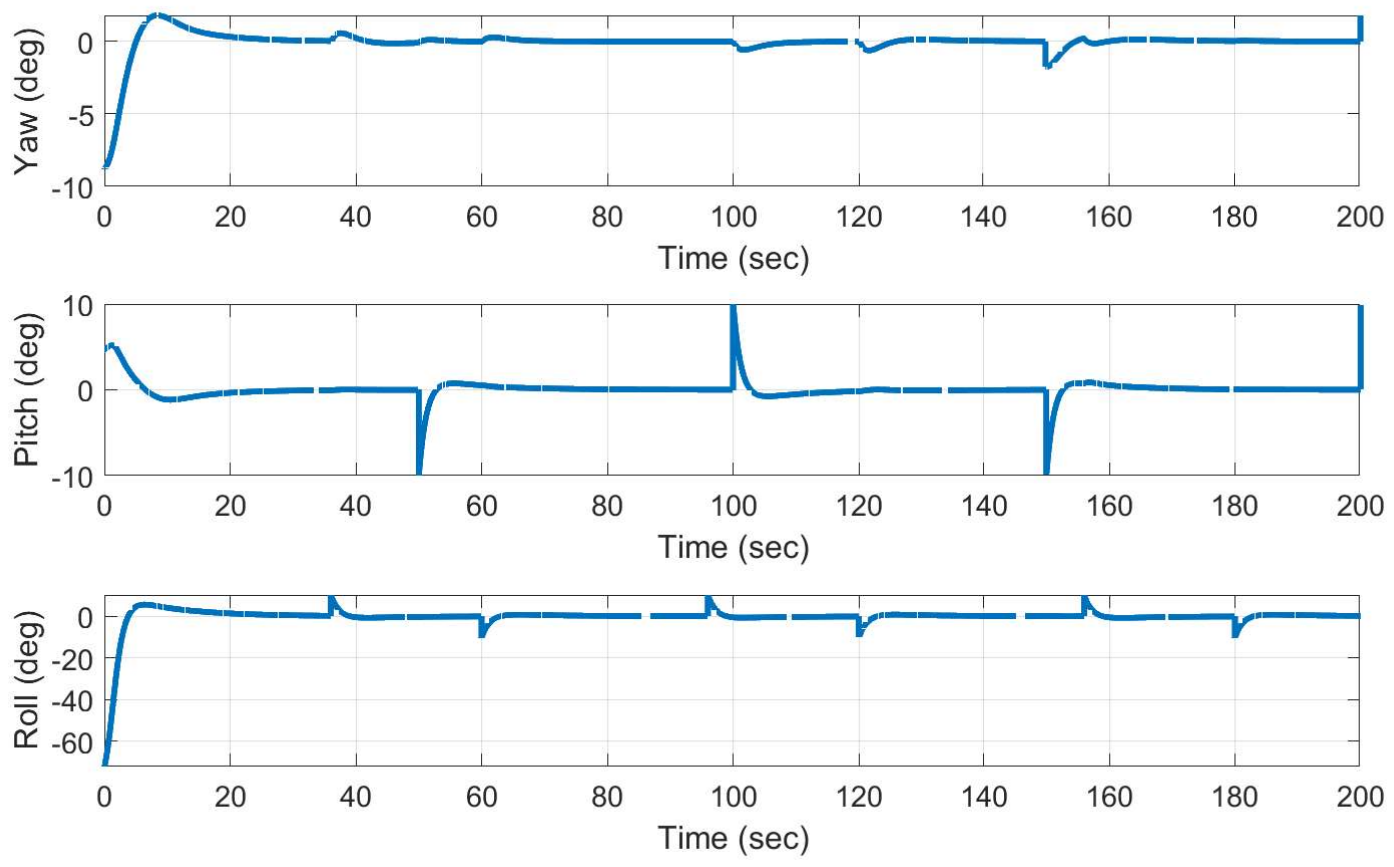}
    \caption{Tracking error for roll, pitch, and yaw angles over time.}
    \label{fig:tracking_error}
\end{figure}

After stabilization, a series of time-varying attitude references were commanded to test the tracking performance of the control strategy. Figure~\ref{fig:tracking_error} presents the attitude tracking error over time for roll, pitch, and yaw. The results demonstrate the effectiveness of the controller in accurately tracking the desired references with negligible steady-state error and minimal transient oscillations.

\begin{figure}[h!]
    \centering
    \includegraphics[width=0.6\textwidth, angle=270]{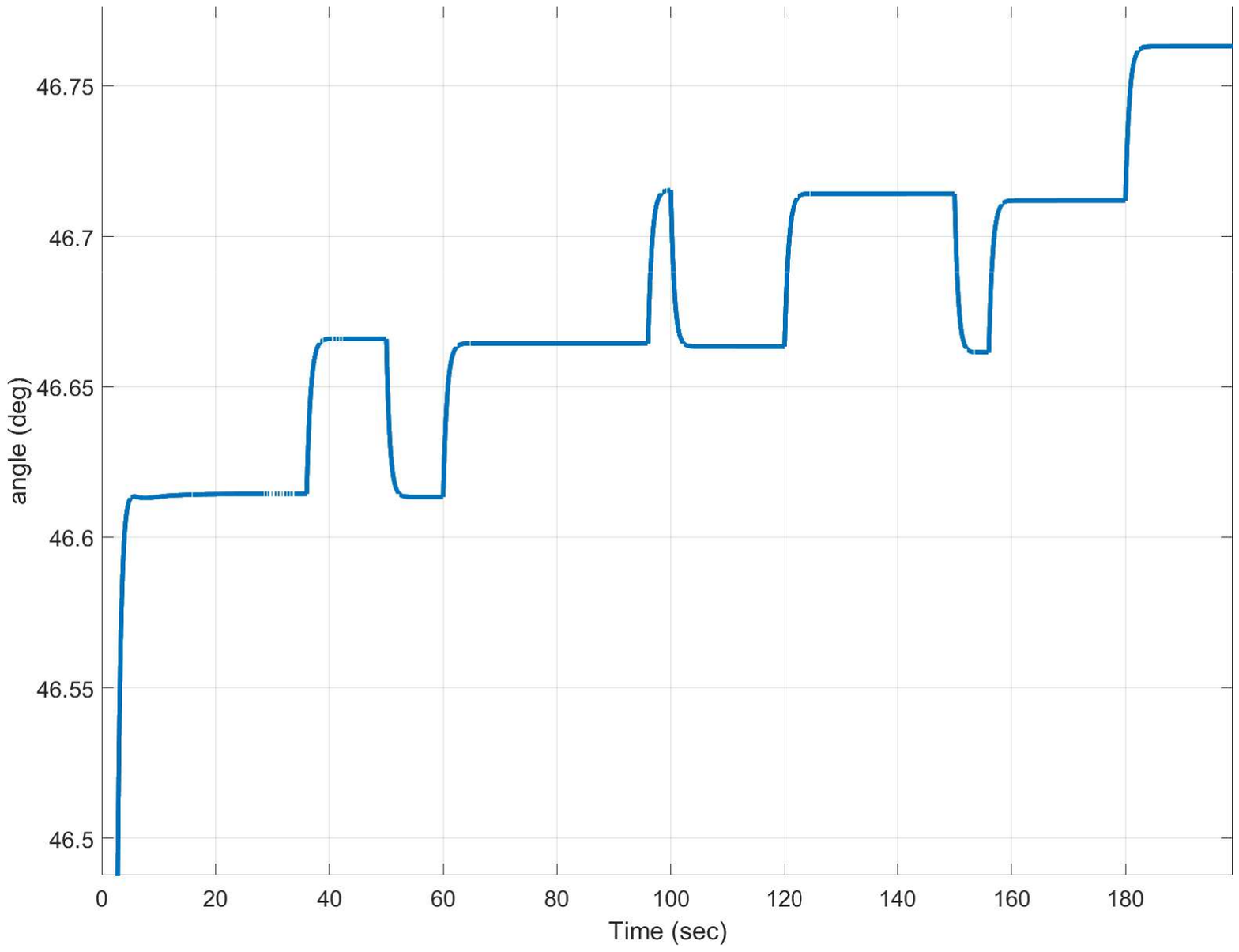}
    \caption{Angle $u$ (internal control - see Figure \ref{drone}) of the quadcopter's arms during the control maneuvers.}
    \label{fig:arm_angles}
\end{figure}

Additionally, the simulation provided a visual representation of the quadrotor behavior during the experiment. The simulated angles of the vehicle’s arms, as depicted in Figure~\ref{fig:arm_angles}, reflect the physical correspondence between the applied control inputs and the resulting system dynamics. This behavior highlights the controller's ability to generate consistent and predictable responses under varying conditions.

\section*{Appendix A: Forced Euler-Poincar\'e and forced Lie-Poisson equations for control-dependent systems}\label{appendixB}

\begin{theorem}
Let $\ell:\mathfrak{g}\times U\to\mathbb{R}$ be a controlled reduced Lagrangian function and $f:G\times\mathfrak{g}\times U\times\bar{U}\to G\times\mathfrak{g}^{*}$ external forces. For a curve $g(t)$ consider the curve on the Lie algebra $\mathfrak{g}$ given by $\omega(t)=g^{-1}\dot{g}$, then Lagrange-d'Alembert principle 
\begin{equation}\label{LDP}
\delta\int_{0}^{T}\ell(\omega(t), u(t))dt=\int_{0}^{T}f(g(t),\omega(t),u(t),\tau(t))\delta g\, dt 
\end{equation} holds for all variations of $g(t)$ constrained to $\delta\omega=\dot{\eta}+ad_{\omega}\eta$, where $\eta(t)$ is a curve on $\mathfrak{g}$, $\eta=g^{-1}\delta g$ with $\eta(0)=\eta(T)=0$, is equivalent to the forced controlled Euler-Poincar\'e equations  \begin{equation}\label{forced}\frac{d}{dt}\left(\frac{\partial \ell}{\partial\omega}\right)-ad_{\omega}^{*}\left(\frac{\partial \ell}{\partial\omega}\right)=T_{e}^{*}\mathcal{L}_gf(g,\omega,u,\tau)\end{equation} together with $\dot{g}(t)=T_{e}\mathcal{L}_g \omega(t)$.
\end{theorem}

\textit{Proof:} Taking variations in \eqref{LDP}, 

\begin{align*}
&\delta\int_{0}^{T}\ell(\omega(t),u(t))\,dt-\int_{0}^{T}f(g(t),\omega(t),u(t),\tau(t))\delta g\, dt\\&=\int_{0}^{T}\Big\langle\frac{\partial \ell}{\partial\omega},\delta\omega\Big\rangle dt-\int_{0}^{T}\Big\langle f,T_{e}\mathcal{L}_g\eta\Big\rangle dt\\
&=\int_{0}^{T}\Big\langle\frac{\partial \ell}{\partial\omega}, \dot{\eta}+ad_{\omega}\eta\Big\rangle-\int_{0}^{T}\Big\langle T_{e}^{*}\mathcal{L}_gf,\eta\Big\rangle\, dt\\
&=\int_{0}^{T}\left(\Big\langle\frac{\partial \ell}{\partial\omega},\dot{\eta}\Big\rangle+\Big\langle\frac{\partial \ell}{\partial\omega},ad_{\omega}\eta\Big\rangle\right)dt-\int_{0}^{T}\Big\langle T_{e}^{*}\mathcal{L}_{g}f,\eta\big\rangle dt\\
&=\int_{0}^{T}\Big\langle-\frac{d}{dt}\left(\frac{\partial \ell}{\partial\omega}\right)+ad_{\omega}^{*}\left(\frac{\partial \ell}{\partial\omega}\right)-T_{e}^{*}\mathcal{L}_{g}f,\eta\Big\rangle dt,
\end{align*} where we have used in the last equality integration by parts and the fact that $\eta(0)=\eta(T)=0$. Therefore the variational principle \eqref{LDP} holds if 
$\frac{d}{dt}\left(\frac{\partial \ell}{\partial\omega}\right)-ad_{\omega}^{*}\left(\frac{\partial \ell}{\partial\omega}\right)=T_{e}^{*}\mathcal{L}_{g}f$.\hfill$\square$

Note that by expanding \eqref{forced}, if the matrix $\left(\frac{\partial^2 \ell}{\partial\omega\partial\omega}\right)$ is non-singular, by the implicit function theorem, we can write \eqref{forced} as an explicit system of first-order ordinary differential equations $\dot{\omega}=K(g,\omega, u, \dot{u},\tau)$ with $K:G\times\mathfrak{g}\times TU\times\bar{U}\to\mathbb{R}$.

If the map $(u,\omega)\mapsto (u,\Pi=\frac{\partial \ell}{\partial\omega})$ is a diffepmorphism we can locally define the controlled reduced Hamiltonian $h:\mathfrak{g}^{*}\times U\to\mathbb{R}$, given by   $$h(\Pi,u)=\langle\Pi,\omega(\Pi)\rangle-\ell(\omega(\Pi),u)$$ with $\omega=\frac{\partial h}{\partial\Pi}$. Therefore, the forced controlled Euler-Poincar\'e equations can be written as the forced controlled Lie-Poisson equations 
\begin{equation}\label{forcedliepoisson}
    \dot{\Pi}=ad^{*}_{\frac{\partial h}{\partial\Pi}}\Pi+T^{*}_{e}\mathcal{L}_{g}f.
\end{equation}

For the Lagrangian of the rigid body with controlled inertia $L:\mathfrak{so}(3)\times U\to\mathbb{R}$ given by $L(\omega,u)=\frac{1}{2}\langle I(u)\omega,\omega\rangle$, the forced Euler-Poincar\'e equations are \begin{equation}\label{eqforL2}
I(u)\dot{\omega}+\left(\frac{\partial I}{\partial u}\dot{u}\right)\omega=F-\omega\times I(u)\omega.
\end{equation} with $F=T_{e}^{*}\mathcal{L}_gf$.

From the Legendre transformation. $\Pi=I\omega$ then $\omega=I^{-1}(u)\Pi$ provided the existence of $I^{-1}(u)$ and hence the reduced hamiltonian function $h:\mathfrak{g}^{*}\times U\to\mathbb{R}$ is given by $h(\Pi,u)=\frac{1}{2}\langle\Pi,I^{-1}(u)\Pi\rangle$. Therefore, the forced controlled Lie-Poisson equations are given by $$\dot{\Pi}=F-\Pi\times I^{-1}(u)\Pi,$$ where now $F:G\times\mathfrak{g}^*\times U\times\bar{U}\to \mathfrak{g}^{*}$.

\subsection*{Appendix B: Constraint algorithm}

In this appendix, we shall introduce the Gotay-Nester-Hinds
algorithm \cite{GoNe,GoNesHinds78} that we have used in the paper. Take the following triple $(M, \Omega, H)$ consisting of a
smooth manifold $M$, a closed 2-form $\Omega$ and a differentiable
function $H: M\rightarrow \mathbb{R}$. On $M$ we consider the equation
\begin{equation}\label{pres}
i_X\Omega=dH
\end{equation}
Since we are not assuming that $\Omega$ is nondegenerate (that is,
$\Omega$ is not, in general, symplectic) then Equation (\ref{pres}) has no solution
in general, or the solutions are not defined everywhere. In the most
favorable case, Equation (\ref{pres}) admits a global (but not
necessarily unique) solution $X$. In this case, we say that the
system admits global dynamics.
 Otherwise, we select the subset of points of $M$, where such a solution exists. We denote by $M_2$ this subset and we will assume that it is a submanifold of  $M=M_1$.
 We are then assured that the equations (\ref{pres}) admit a solution $X$ defined at all points of $M_2$, but $X$ need not be tangent to $M_2$, hence, does not necessarily induce a dynamics on $M_2$.
  So we impose an additional tangency
condition, and we obtain a new submanifold $M_3$ along which there
exists a solution $X$, but, however, such an $X$ need not be tangent
to $M_3$. Continuing this process, we obtain a sequence of
submanifolds
\[
\cdots M_s \hookrightarrow \cdots \hookrightarrow M_2 \hookrightarrow M_1=M
\]
where the general description of $M_{l+1}$ is
\[
\begin{array}{rl}
M_{l+1}:=\{p\in M_{l} \textrm{ such that there exists } X_p\in T_pM_l\,
\textrm{ satisfying }
 i_{X_p}\Omega(p)=dH(p) \}.
\end{array}
\]
If the algorithm terminates at a nonempty set, in the sense that  at
some $s\geq 1$ we have  $M_{s+1}=M_s$. Then we say that $M_s$ is the
final constraint submanifold which is denoted by $M_f$. It may still
happen that $\dim M_f=0$, that is, $M_f$ is a discrete set of
points, and in this case the system does not admit a proper
dynamics. But in the case when $\dim M_f>0$, there exists a
well-defined solution $X$ of (\ref{pres}) along $M_f$.



\bibliographystyle{siamplain}
\bibliography{references}
\end{document}